\documentclass{article}

\usepackage[preprint,nonatbib]{neurips_2023}

\usepackage[utf8]{inputenc} 
\usepackage[T1]{fontenc}    
\usepackage{hyperref}       
\usepackage{url}            
\usepackage{booktabs}       
\usepackage{amsfonts}       
\usepackage{nicefrac}       
\usepackage{microtype}      
\usepackage{xcolor}         

\usepackage{amsmath}
\usepackage{amssymb}
\usepackage{mathtools}
\usepackage{amsthm}
\usepackage{relsize}
\usepackage{wrapfig}
\usepackage{mathrsfs}
\usepackage{yfonts}
\usepackage[capitalize,noabbrev]{cleveref}

\usepackage{tikz-network}
\usepackage[figurename=Fig.]{caption}
\usepackage{subcaption}

\usepackage{graphicx}
\usepackage[english]{babel}
\usepackage{ amsfonts}
\usepackage{amssymb}
\usepackage{mathtools, url}
\usepackage{graphicx, balance}
\usepackage{changepage}
\usepackage{bm}
\usepackage{lineno}
\usepackage{algorithm,algorithmic}

\usepackage{enumitem}

\usepackage{xcolor}
\hypersetup{
    colorlinks,
    linkcolor={red},
    citecolor={green},
    urlcolor={blue!80!black}
}

\addto\captionsenglish{}
\theoremstyle{plain}
\newtheorem{theorem}{Theorem}[section]

\newtheorem{lemma}[theorem]{Lemma}

\theoremstyle{definition}
\newtheorem{definition}[theorem]{Definition}

\theoremstyle{remark}

\newcommand\bigforall{\mbox{\Large $\mathsurround=0pt\forall$}}

\title{Uncertainty Informed Optimal Resource Allocation with Gaussian Process based Bayesian Inference
}

\author{%
  Samarth Gupta \\
  Massachusetts Institute of Technology, USA\\
  \texttt{samarthg@mit.edu} \\
  \And
   Saurabh Amin \\
   Massachusetts Institute of Technology, USA\\
   \texttt{amins@mit.edu} \\
}


\begin{document}
\setlength{\abovedisplayskip}{1pt}
\setlength{\belowdisplayskip}{1pt}

\maketitle

\begin{abstract}

We focus on the problem of uncertainty informed  allocation of medical resources (vaccines) to heterogeneous populations for managing epidemic spread. We tackle two related questions: (1) For a compartmental ordinary differential equation (ODE) model of epidemic spread, how can we estimate and integrate parameter uncertainty into resource allocation decisions? (2) How can we computationally handle both nonlinear ODE constraints and parameter uncertainties for a generic stochastic optimization problem for resource allocation? To the best of our knowledge current literature does not fully resolve these questions. Here, we develop a data-driven approach to represent parameter uncertainty accurately and tractably in a novel stochastic optimization problem formulation. We first generate a tractable scenario set by estimating the distribution on ODE model parameters using Bayesian inference with Gaussian processes. Next, we develop a parallelized solution algorithm that accounts for scenario-dependent nonlinear ODE constraints. Our scenario-set generation procedure and solution approach are flexible in that
they can handle any compartmental epidemiological ODE model. 
Our computational experiments on two different non-linear ODE models (SEIR and SEPIHR) indicate that accounting for uncertainty in key epidemiological parameters can improve the efficacy of time-critical allocation decisions by 4-8\%. This improvement can be attributed to data-driven and  optimal (strategic) nature of vaccine allocations, especially in the early stages of the epidemic when the allocation strategy can crucially impact the long-term trajectory of the disease.

\end{abstract}

\section{Introduction}
\vspace{-0.3cm}

In this paper we study the problem of \emph{uncertainty informed} optimal resource allocation to control the spread of an infectious disease such as Covid-19.
We develop a \emph{data-driven, scalable} and ODE \emph{model agnostic} approach while \emph{accounting for uncertainty} for the vaccine allocation problem.
Our approach is flexible in that it can be easily adapted to other control strategies such as imposing lock-downs \cite{candogan} and allocation of other resources such as medical personnel, supplies, testing facilities \& etc.

The vaccine allocation problem  has been well studied in the literature. This includes earlier works like \cite{model1967,becker1975} to more recent optimization based methods like \cite{bertsimas_vaccine,chenyi_melvyn}.
Researchers have also studied ways 
to incorporate uncertainty through stochastic epidemiological modelling \cite{clancy_2007, chenyi_melvyn} 
or stochastic optimization with uncertain parameters \cite{tanner_vacc_stochastic, vacc_2phase_stochastic}.
 However, prior works 
 have two major limitations:\begin{enumerate}[leftmargin=0.65cm]
    \item Most papers such as \cite{tanner_vacc_stochastic, vacc_2phase_stochastic, esra_multi_stage_stoc,clancy_2007} which claim to account for uncertainty, do not provide a principled \emph{data-driven} method to model (and estimate)  uncertainty. They simply model the allocation problem as a stochastic program under the assumption that a scenario-set exists without outlining a principled procedure on how to generate or estimate this scenario-set from \emph{data}. 
    Clearly, this does not effectively solve the problem of uncertainty informed vaccine allocation.
    
    \item The presence of product term between the susceptible (S) and infected (I) population is a key characteristic of most compartmentalized epidemiological ODE models \cite{bertsimas_vaccine}. Due to this non-linearity, the resource allocation problem with the discretized ODEs results in a non-convex quadratic program. This is difficult to solve even in the nominal case i.e. without accounting for uncertainty, let alone uncertainty informed. To avoid the product term, previous papers \cite{tanner_vacc_stochastic, vacc_2phase_stochastic, esra_multi_stage_stoc,esra_ventilator} resort to using simple (linear) epidemiological models so that the discretized ODEs result in a linear program which is easy to solve. 
    Such linear models are limited in their ability to capture the true underlying non-linear dynamics of disease transmission;
    hence the resulting allocation strategies are not globally optimal. 
\end{enumerate}

In this work, we address both of the above limitations by making following novel contributions: 
\vspace{-0.2cm}
\begin{enumerate}[leftmargin=0.65cm]
    \item We make progress in resolving the issue of incorporating parameter uncertainty in the resource (vaccine) allocation  problem in a \emph{data-driven manner}. We do this by making connections with the  ODE parameter estimation literature with Bayesian inference using GPs with gradient matching methods. We show that the posterior-distributions can be used to represent uncertainty through a tractable scenario-set by formulating the scenario-reduction problem as an optimal-transport problem. This optimal transport problem can be easily solved using k-means clustering \cite{rujeerapaiboon_sc_reduction}.  

    \item We provide a novel formulation for the uncertainty informed vaccine allocation problem as a stochastic optimization problem. We develop technical results for the feasibility and decomposability of this stochastic program. Using these technical results, we develop an \emph{parallelized}, scalable iterative solution algorithm to solve the stochastic program while retaining the original \emph{non-linear}, continuous-time  ODE model constraints. 
Due to this ODE \emph{model agnostic} nature of our approach, we are also able to account for different levels of mobility within different sub-populations and the temporal variations in the onset of the pandemic in each of these sub-populations.

   \item We provide extensive empirical results on two different ODE models (i.e. the SEIR and the SEPIHR models) in sections \ref{sec:empirical_sampling_result}, \ref{sec:sce_reduction}, \ref{sec:exp_simu_results} and Supplementary Information (SI). Our results demonstrate that with optimal vaccine allocation, peak infections can be reduced by around 35\%. More importantly, a further gain of around \textbf{4} to \textbf{8}\% can be achieved when incorporating uncertainty. 
\end{enumerate}

\vspace{-0.3cm}
\section{Epidemiological Modelling and Pitfalls of Classical Parameter estimation} \label{sec:model}
\vspace{-0.25cm}
Mathematical modelling of pandemics (including epidemics) has an extensive literature going back to 1960s \cite{model1967, becker1975}. A fairly recent and concise overview can be found at \cite{brauer_overview}. Throughout literature, modelling the spread of different diseases using a compartmentalized model through a set of time-dependent ordinary differential equations (ODEs) is  common and widely used \cite{compartmental_models}. 
Following the recent literature on covid-19 \cite{delphi, acemoglu_parise, parise_repro_no,pnas_ensemble} we also adopt the compartmentalized modelling approach. 

\begin{figure}[h!]
\vspace{-0.2cm}
\begin{minipage}[t]{0.25\linewidth}
       \centering 
    \includegraphics[scale = 0.40]{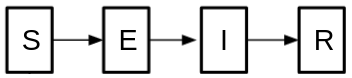}
     \vspace{-0.5cm}
     \caption{ SEIR model.}\label{fig:seir}
\end{minipage}
\begin{minipage}{0.74\linewidth}
\begin{small}
    \begin{equation}\label{eq:seir}
    \begin{rcases}
    \begin{alignedat}{10}
    &\frac{dS(t)}{dt}   \coloneqq \dot{S}(t) =  - \frac{\alpha}{N} S(t) I(t), \frac{dI(t)}{dt}  \coloneqq \dot{I}(t) = \beta E(t) - \gamma  I(t)\\        
    &\frac{dE(t)}{dt} \coloneqq \dot{E}(t) = \frac{\alpha}{N} S(t) I(t) -  \beta E(t), \frac{dR(t)}{dt}   \coloneqq \dot{R}(t) = \gamma I(t) 
    \end{alignedat}
    \end{rcases}\hspace{-0.25cm}
    \end{equation}        
\end{small}    
\end{minipage}
\vspace{-0.3cm}
\end{figure}

A popular epidemiological model which we use is the SEIR model, shown in fig. \ref{fig:seir}. In this model the entire population (of size N) is divided into
four states: Susceptible (S), Exposed (E), Infected (I) and Recovered (R). The evolution of each state or the system dynamics is governed by equations in \eqref{eq:seir}.

In \eqref{eq:seir}, $\alpha, \beta, \text{ and } \gamma$ are the model parameters and control the rate at which fraction of the population moves from one compartment to another. These  model parameters are to be  estimated from the available time-series data which we discuss subsequently.
Mobility levels can be easily  incorporated by adjusting the infection rate $\alpha$ accordingly. 


\begin{wrapfigure}{r}{0.27\linewidth}
\vspace{-0.5cm}
\captionsetup{singlelinecheck=on, margin={0.0cm, 0.0cm}, justification=justified, format=plain}
\centering
\resizebox{\linewidth}{!}{%
    \includegraphics[angle=-0,origin=c, scale = 0.4]{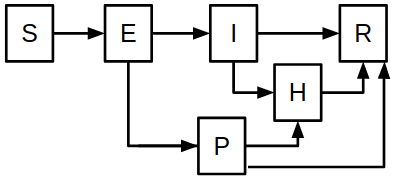}}
    \vspace{-0.5cm}
    \centering
    \caption{\small{SEPIHR model }}\label{fig:sepihrm_model}    
\vspace{-0.4cm}
\end{wrapfigure}

Note that we use SEIR model only as a prototypical model, however, all our subsequent discussion including technical results and solution algorithm holds true for other ODE based models as well. In fact, in addition to the SEIR model, we also provide results on a second model, i.e. the SEPIHR model \cite{sepihrm} with additional states P (for protective quarantine) and H (for hospitalised quarantined) shown in fig. \ref{fig:sepihrm_model}. The functional form of the ODEs for this model are provided in SI.  

Given the time-series data such as number of daily infections and deaths, the main question arises how to estimate SEIR model parameters i.e.  $\alpha, \beta, \text{ and } \gamma$  from this data. Therefore, we next discuss the commonly used non-linear least squares approach for ODE parameter estimation and its associated pitfalls, thus providing motivation for adopting Bayesian viewpoint.

\vspace{-0.25cm}
\subsection{Classical Parameter Estimation: Non-linear Least Squares (NLLS)} \label{sec:param_estimation_nlls}
\vspace{-0.25cm}



Before describing the NLLS approach, we briefly describe the \emph{initial-value problem} (IVP) in the context of ODEs.  For a given (or fixed) set of parameter values and initial conditions (denoted $\textbf{x}_0$), a systems of ODEs can be numerically solved using an off-the-shelf ODE solver such as ODE45 in matlab or ODEINT in python. 
The solved  system (also referred to as simulation) provides the value (or estimates) of different states  at the specified time-stamps. 

For a given set of parameters, using the estimated state values obtained by solving the IVP and time-series data, discrepancy or the least-squares error can be computed. This can be turned into a optimization problem where we want to find those values of the model parameters for which the least-squares error is minimized. L-BFGS is commonly used to solve such problems \cite{delphi}.
Mathematically for SEIR model,  the NLLS problem can be written as follows:
\begin{small}
\begin{subequations}
\begin{alignat}{10}
\displaystyle \min_{{\alpha, \beta, \gamma}}  \sum_{t=1}^{N} & \bigg( (y_R^t - R(t) )^2  +  (y_I^t - I(t))^2 \bigg)  \notag \\
\text{ s.t. } & \{ \eqref{eq:seir}  \} \;\;\; \forall \; t \in \{1,\dots,N \}  \text{ and } [S(0) , E(0),  I(0), R(0)] = \textbf{x}_0 \notag
\end{alignat}
\end{subequations}
\end{small}
where $y_R^t \text{ and } y_I^t$ denote count data for infected and removed individuals at time $t$. The optimal parameters obtained after solving NLLS can then be used to re-solve the ODE system to make predictions for future time as well.   
\vspace{-0.15cm}
\paragraph{Why account for Uncertainty? }
NLLS discussed above can provide sufficiently reliable \emph{point estimates} of the parameter values and  predictions of new cases into the future provided the \emph{time-series data is accurate}. Using these point estimates resource (vaccine) allocation problem is to be solved subsequently.
The efficacy of the overall allocation solution in real-world is highly dependent on the accuracy of the predicted point estimates which are only as good as the data from which these estimates are generated. For Covid-19, the data reported by various private organizations and government agencies can be severely biased, under-reported \cite{covid_deaths_bias} and erroneous due to numerous reasons \cite{jordan_bias}. Reliance on these point estimates can result in severe region-wide inefficiencies. To address these issues and also account for potential modelling errors, we  incorporate \emph{uncertainty} through Bayesian inference to estimate the \emph{joint-distribution} of ODE model parameters from data, which we discuss next.

\vspace{-0.35cm}
\section{Bayesian Parameter Estimation}\label{sec: param_estimation_bayesian}
\vspace{-0.3cm}

Bayesian inference for estimating ODE parameters has been well studied in literature \cite{ramsay}, however in the absence of closed form posterior and
the requirement of solving the ODE system in each sampling iteration  makes inference difficult. To overcome this limitation, \cite{calderhead} proposed the use of Gaussian Processes (GPs) to model the evolution of a state over time while exploiting the fact that derivative of a Gaussian process is also a Gaussian process. This significantly helps in achieving tractability and allows Bayesian inference to be computationally feasible. Following \cite{calderhead}, numerous other related works like \cite{dondelinger, barber, controversy, niu, nico, fgpgm, odin} have been proposed which also employ the use of GPs to efficiently estimate the 
parameters of a non-linear ODE system (for eg: SEIR model). We discuss some of these works, in particular the approach of \cite{fgpgm} which is useful to our problem setting.

Consider a set of $K$ time-dependent states denoted as $\textbf{x}(t)= [x_1(t),\dots,x_K(t)]^T$.
The evolution of each of these $K$ state over time is defined by a set of $K$ time-dependent arbitrary differential equations denoted as follows:
\begin{equation}\label{eq:def_ode}
    \displaystyle \dot{x}_i(t) = \frac{d{x}_i(t)}{dt} =  \text{f}_i(x(t),\theta,t)   \quad \forall \quad i \in \{ 1, \dots, K\}
\end{equation}

where the functional form of $\text{f}_i$ is known (for eg. SEIR model). Noisy observations (i.e. the time-series data) of  each of the $K$ states (denoted $\textbf{y}(t) = [ y_1(t),\dots,y_K(t)  ]^T$) at $N$ different time points where $t_1 < \dots < t_N $ are available, i.e. 
\begin{equation*}
\begin{rcases}
\begin{matrix}
    y_1(t)  =  x_1(t)  +  \epsilon_1(t) \\ 
    \vdots   \\ 
    y_K(t)  =  x_K(t) +  \epsilon_K(t) \\ 
\end{matrix} 
\end{rcases} \;, 
\begin{aligned}
  &\text{ where } \epsilon_i(t) \sim \mathcal{N}(0,\sigma^2_i). \\  
  & \;\; \forall \;\; t \in  \{t_1 ,\dots,t_N\}
\end{aligned}  
\end{equation*}
Let $\bm{\epsilon}(t) = [ \epsilon_1(t),\dots,\epsilon_K(t)  ]^T$, then in vector notation we have $\textbf{y}(t) = \textbf{x}(t) + \bm{\epsilon}(t) $.
As there are  $N$ observations for each of the $K$ states, for a clear exposition we introduce matrices of size $K \times N$ as follows: $\textbf{X}  = [ \textbf{x}(t_1), \dots, \textbf{x}(t_N) ]$ and $\textbf{Y}  = [ \textbf{y}(t_1), \dots, \textbf{y}(t_N) ]$.
We can then write :
\begin{alignat}{10}
\displaystyle P(\textbf{Y}| \textbf{X}, \sigma) &= \prod_{k} \prod_{t} P(y_k(t)| x_k(t),\sigma) =\prod_{k} \prod_{t} \mathcal{N}(y_k(t)| x_k(t),\sigma^2) \label{eq:y_given_x_sigma}
\end{alignat}
 
\cite{calderhead} proposed placing a Gaussian process prior on $\textbf{x}_k$. 
Let $\bm{\mu}_k$ and $\bm{\phi}_k$ be the hyper-parameters of this Gaussian process, we can then write:
\begin{equation}\label{eq:gp}
p(\textbf{x}_k|\bm{\mu}_k,\bm{\phi}_k) = \mathcal{N}(\textbf{x}_k | \bm{\mu}_k, \textbf{C}_{\phi_{k}})
\end{equation}

In \eqref{eq:gp}, $\textbf{C}_{\phi_{k}}$, denotes the Kernel (or the covariance) matrix for a predefined kernel function with hyper-parameters $\bm{\phi}_k$. 
As differentiation is a linear operator therefore the derivative of a Gaussian process is also a Gaussian process (see ch-9 in \cite{rasmussen}  and \cite{solak}). Therefore a Gaussian process is closed under differentiation and the joint distribution of the state variables $\textbf{x}_k$ and their derivatives $\dot{\textbf{x}}_k$ is a multi-variate Gaussian distribution as follows:
\begin{equation}\label{eq:joint_x_xdot}
\begin{bmatrix}
\textbf{x}_k \\
\dot{\textbf{x}}_k
\end{bmatrix} \sim \mathcal{N}\bigg(
\begin{bmatrix}
\bm{\mu}_k \\
\bm{0}
\end{bmatrix} ,
\begin{bmatrix}
 \textbf{C}_{\phi_{k}} , '\textbf{C}_{\phi_{k}} \\
 \textbf{C}'_{\phi_{k}}, \textbf{C}''_{\phi_{k}}
\end{bmatrix}
\bigg)
\end{equation}

where $\textbf{C}_{\phi_{k}}$ and $\textbf{C}''_{\phi_{k}}$ are the kernel matrices for the state $\textbf{x}_k$ and its derivative $\dot{\textbf{x}}_k$ respectively, while  $'\textbf{C}_{\phi_{k}} \text{ and } \textbf{C}'_{\phi_{k}}$ are the cross-covariance kernel matrices between the states and their derivatives. The functional form of the entries of $\textbf{C}_{\phi_{k}}, \textbf{C}''_{\phi_{k}}, 
'\textbf{C}_{\phi_{k}} \text{ and } \textbf{C}'_{\phi_{k}} $ are provided in SI. Importantly, this implies that using the Gaussian process defined on the state variables $\textbf{x}_k$, we can also make predictions about their derivatives  $\dot{\textbf{x}}_k$. 
From \eqref{eq:joint_x_xdot}, we can compute the conditional distribution of the state derivatives as:
\begin{equation}\label{eq:graph_2}
p(\dot{\textbf{x}}_k | \textbf{x}_k, \bm{\mu}_k, \bm{\phi}_k)  = \mathcal{N}(\dot{\textbf{x}}_k |\textbf{m}_k,\textbf{A}_k )
\end{equation}
where $\textbf{m}_k = '\textbf{C}_{\phi_{k}}{\textbf{C}_{\phi_{k}}}^{-1} (\textbf{x}_k - \bm{\mu}_k ); \;\textbf{A}_k = \textbf{C}''_{\phi_{k}} - '\textbf{C}_{\phi_{k}}{\textbf{C}_{\phi_{k}}}^{-1}\textbf{C}'_{\phi_{k}}$. Note that $p(\dot{\textbf{x}}_k | \textbf{x}_k, \bm{\mu}_k, \bm{\phi}_k)$ corresponds to the second, i.e. GP part of the graphical model in fig. \ref{fig:ode_gp_model}.

\begin{wrapfigure}{r}{0.35\linewidth}
\vspace{-0.3cm}
    \SetVertexStyle[TextFont= \Large]
    \centering
    \resizebox{\linewidth}{!}{%
    \begin{tikzpicture}
    
    \begin{pgfonlayer}{background}
     \draw[color =red, style= dashed, rounded corners=0.2cm] (-0.75, -0.6) rectangle (3.3, 2.1) {};
     \draw[color =blue,style= dashed, rounded corners=0.2cm] (3.5, -0.6) rectangle (7.5, 2.1) {};
    \end{pgfonlayer}
    
    \Vertex[ size= 0.8,opacity =0, label = $\dot{\textbf{x}}$,   x=1.3, y=0 ]{x_dot_ode}
    \Vertex[ size= 0.8,opacity =0, label = $\bm{\theta}$,        x=-0.2,   y=0 ]{theta}
    \Vertex[ size= 0.8,opacity =0, label = $\textbf{x}$,         x=2.8,   y=0]{x_ode}
    \Vertex[ size= 0.8,opacity =0, label = $\lambda$,             x=1.3, y=1.5 ]{lambda}
    \Edge[Direct=true,color= black, lw =1pt](theta)(x_dot_ode)
    \Edge[Direct=true,color= black,lw =1pt](lambda)(x_dot_ode)
    \Edge[Direct=true,color= black,lw =1pt](x_ode)(x_dot_ode)
    \node at (1.3,-1) { {\Large ODE model}};
    
    \Vertex[ size= 0.8,opacity =0, label = $\dot{\textbf{x}}$,   x=4, y=1.5 ]{x_dot}
    \Vertex[ size= 0.8,opacity =0, label = $\textbf{x}$,         x=5.5, y=1.5]{x}
    \Vertex[ size= 0.8,opacity =0, label = $\textbf{y}$,             x=7, y=1.5 ]{y}
    \Vertex[ size= 0.8,opacity =0, label = $\phi$,               x=5.5, y=0 ]{phi}
    \Vertex[ size= 0.8,opacity =0, label = $\sigma$,             x=7, y=0 ]{sigma}
    \Edge[Direct=true,color= black, lw =1pt](phi)(x)
    \Edge[Direct=true,color= black,lw =1pt](phi)(x_dot)
    \Edge[Direct=true,color= black,lw =1pt](sigma)(y)
    \Edge[Direct=true,color= black,lw =1pt](x)(y)
    \Edge[Direct=true,color= black,lw =1pt](x)(x_dot)
    \node at (5.5,-1) { {\Large GP model}};
    
    \Edge[Direct=false,style =dashed, color= black, lw =1pt](x_dot_ode)(x_dot)
    \Edge[Direct=false,style=dashed, color= black,lw =1pt](x_ode)(x)
    
    \end{tikzpicture}
    }
      \vspace{-0.5cm}
        \caption{} \label{fig:ode_gp_model}
        \vspace{-0.4cm}
\end{wrapfigure}

Using the functional form of the ODE system in \eqref{eq:def_ode} and with state specific Gaussian additive noise $\lambda_k$, we can write
\begin{equation}\label{eq:graph_1}
    p(\dot{\textbf{x}}_k| \textbf{X},\bm{\theta},\lambda_k) = \mathcal{N}( \dot{\textbf{x}}_k |\textbf{f}_k(\textbf{X}, \bm{\theta}) , \lambda_k\textbf{I})
\end{equation}

where $\textbf{f}_k(\textbf{X}, \bm{\theta}) = [ \text{f}_k(\text{x}(t_1), \bm{\theta} ),\dots,\text{f}_k(\text{x}(t_N), \bm{\theta} )   ]^T$.
Note that \eqref{eq:graph_1} corresponds to the ODE part of the graphical model in the fig. \ref{fig:ode_gp_model}.

\SetVertexStyle[TextFont= \Large]
\begin{wrapfigure}{r}{0.25\linewidth}
\vspace{-0.35cm}
\centering
  \captionsetup{singlelinecheck=on, margin={0.0cm, 0.0cm}, justification=justified, format=plain}
\resizebox{0.9\linewidth}{!}{%
\begin{tikzpicture}
\Vertex[ size= 0.8,opacity =0, label = $\dot{\textbf{x}}$,   x=4, y=1.5 ]{x_dot}
\Vertex[ size= 0.8,opacity =0, label = $\textbf{x}$,         x=5.5, y=1.5]{x}
\Vertex[ size= 0.8,opacity =0, label = $\textbf{y}$,             x=7, y=1.5 ]{y}
\Vertex[ size= 0.8,opacity =0, label = $\phi$,               x=5.5, y=0 ]{phi}
\Vertex[ size= 0.8,opacity =0, label = $\sigma$,             x=7, y=0 ]{sigma}
\Vertex[ size= 0.8,opacity =0, label = $\text{F}_1$,   x=5.5, y=3 ]{F_1}
\Vertex[ size= 0.8,opacity =0, label = $\text{F}_2$,   x=4, y=3 ]{F_2}
\Vertex[ size= 0.8,opacity =0, label = $\lambda$,             x=2.5, y=3 ]{lambda}
\Vertex[ size= 0.8,opacity =0, label = $\bm{\theta}$,        x=7,  y=3 ]{theta}

\Edge[Direct=true,color= black, lw =1pt](phi)(x)
\Edge[Direct=true,color= black,lw =1pt](phi)(x_dot)
\Edge[Direct=false,color= black,lw =1pt](F_1)(F_2)
\Edge[Direct=true,color= black,lw =1pt](lambda)(F_2)
\Edge[Direct=true,color= black,lw =1pt](theta)(F_1)
\Edge[Direct=true,color= black,lw =1pt](x_dot)(F_2)
\Edge[Direct=true,color= black,lw =1pt](x)(F_1)
\Edge[Direct=true,color= black,lw =1pt](sigma)(y)
\Edge[Direct=true,color= black,lw =1pt](x)(y)
\Edge[Direct=true,color= black,lw =1pt](x)(x_dot)
\end{tikzpicture}
}
  \vspace{-0.2cm}
  \begin{small}
    \caption{Combined model.}
    \label{fig:fgpgm_model}
  \end{small}  
\end{wrapfigure}
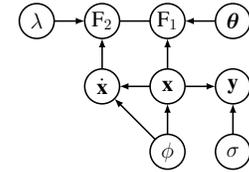
The two models $p(\dot{\textbf{x}}_k | \textbf{x}_k, \bm{\mu}_k, \bm{\phi}_k)$ in \eqref{eq:graph_2} and $p(\dot{\textbf{x}}_k| \textbf{X},\bm{\theta},\lambda_k)$ in \eqref{eq:graph_1}
are combined through 
two new random variables $\textbf{F}_1$ and $\textbf{F}_2$,
resulting in the graphical model shown in fig. \ref{fig:fgpgm_model} \cite{fgpgm}.
Considering a single state (for notational simplicity), for given values of $\textbf{x}$ and $\bm{\theta}$, $\textbf{F}_1$ in fig. \ref{fig:fgpgm_model}, represents the deterministic output of the ODEs, i.e. $\textbf{F}_1= \textbf{f}(\bm{\theta}, \textbf{x})$. The value of $p(\textbf{F}_1| \bm{\theta}, \textbf{x})$ can be written using the Dirac-delta function (denoted $\delta(\cdot)$) as following:
\begin{equation}
p(\textbf{F}_1| \bm{\theta}, \textbf{x}) = \delta(\textbf{F}_1- \textbf{f}(\bm{\theta}, \textbf{x})) 
\end{equation}

Under the assumption that the GP model would be able to capture both, the true states and their derivatives perfectly, then it would imply that $\textbf{F}_1$ is same as $\dot{\textbf{x}}$, i.e.  $\textbf{F}_1 = \dot{\textbf{x}}$. But clearly this assumption is unlikely to hold, therefore to account for any possible mismatch and small error in  the GP states and GP derivatives, this condition is relaxed so that:
\begin{equation}\label{eq:equality_f1_f2}
 \textbf{F}_1 = \dot{\textbf{x}} + \epsilon \eqqcolon \textbf{F}_2, \quad \text{ where } \quad \epsilon \sim \mathcal{N}(\bm{0}, \lambda \textbf{I})  
\end{equation}

The above argument regarding the the error in the states and derivatives of the GP model is captured in the graphical model (fig. \ref{fig:fgpgm_model}) through the use of the random variable $\textbf{F}_2$. From a given state-derivative $\dot{\textbf{x}}$ obtained from the GP model, $\textbf{F}_2$ is obtained after addition of Gaussian noise with standard deviation $\lambda$. The probability density of $\textbf{F}_2$ can then be written as    
\begin{equation}
p(\textbf{F}_2|\dot{\textbf{x}}, \lambda) = \mathcal{N}(\textbf{F}_2 | \dot{\textbf{x}}, \lambda \textbf{I})) 
\end{equation}
Note that the equality constraint in \eqref{eq:equality_f1_f2} is encoded in the graphical model using an un-directed edge between $\textbf{F}_1$ and $\textbf{F}_2$. For the purpose of  inference, this equality constraint is incorporated in the joint density via the Dirac-delta function, i.e. $\delta(\textbf{F}_1 - \textbf{F}_2)$. The joint-density of the whole graphical model (fig. \ref{fig:fgpgm_model}) is given as:
\begin{alignat}{10}
p(\textbf{x}, \dot{\textbf{x}}, \textbf{y}, \textbf{F}_1,\textbf{F}_2, \bm{\theta}| \phi, \sigma, \lambda ) = & p(\bm{\theta}) p(\textbf{x}|\phi) p(\dot{\textbf{x}} | \textbf{x},\phi) p(\textbf{y}| \textbf{x}, \sigma) p(\textbf{F}_1| \bm{\theta}, \textbf{x}) & p(\textbf{F}_2| \dot{\textbf{x}}, \lambda \textbf{I})\delta( \textbf{F}_1-\textbf{F}_2)
\end{alignat}

Finally, the marginal distribution of $\textbf{x},\bm{\theta}$ takes the following form:
\begin{alignat}{10}
&p(\textbf{x},\bm{\theta} |\textbf{y}, \bm{\phi}, \sigma, \lambda ) =  p(\bm{\theta}) \times \mathcal{N}(\textbf{x}| \bm{\mu}, \textbf{C}_{\phi}) \times \mathcal{N}(\textbf{y}|\textbf{x}, \sigma^2 \textbf{I} )\times \mathcal{N}(\textbf{f}(\textbf{x}, \bm{\theta})|\textbf{m}, \textbf{A} + \lambda \textbf{I} ) \label{eq:fgpgm}
\end{alignat}



\vspace{-0.25cm}
\subsection{Empirical Sampling Results} \label{sec:empirical_sampling_result}
\vspace{-0.25cm}
We now provide sampling results on the two disease-transmission ODE models: 1) SEIR model (fig. \ref{fig:seir})  and 2) SEPIHR model (fig. \ref{fig:sepihrm_model}).

\begin{wrapfigure}{r}{0.65\linewidth}
\vspace{-0.1cm}
  \centering
    \captionsetup{singlelinecheck=on, margin={0.25cm, 0.5cm}, justification=justified, format=plain}
  \includegraphics[width=\linewidth]{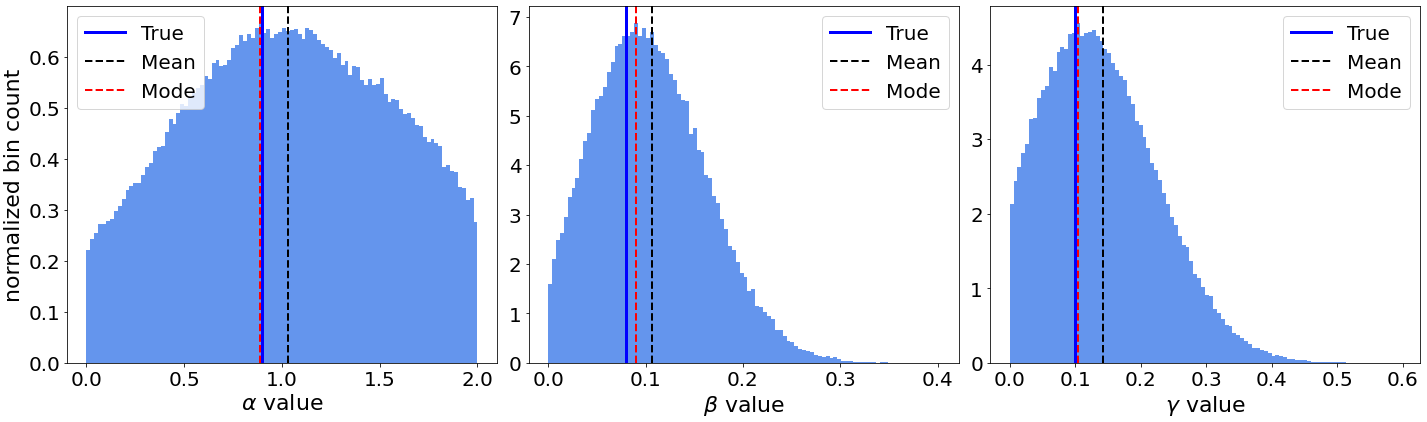}
  \vspace{-0.45cm}
  \captionof{figure}{SEIR: Empirical distribution after $3\times 10^5$ samples from the MCMC sampling procedure.}
  \label{fig:mcmc_dist}    
  \vspace{-0.2cm}
\end{wrapfigure}

\textbf{SEIR:} Using $\alpha = 0.9, \beta = 0.08 \text{ and } \gamma = 0.1$ (as true parameter values) we simulate the SEIR model (eq:\eqref{eq:seir}) to get state values. We add zero-mean Gaussian noise with $\sigma = 0.1$ to each of the simulated state values to generate our dataset. Using the data only for first 15 days ($T=\{1,\dots,15\}$),  we estimate the GP hyper-parameters for states using maximum-likelihood (see SI for details). 
We then run the Metropolis-Hastings MCMC sampling procedure using the density from eq: \eqref{eq:fgpgm} to get our empirical posterior joint-distribution on $\alpha, \beta \text{ and } \gamma$. After removing the burn-in samples, for the remaining  $3\times 10^5$ samples, we plot the marginal distributions along with their mean and mode in fig. \ref{fig:mcmc_dist}. 

\textbf{SEPIHR:} This model has 5 parameters, i.e. $\alpha, \beta, \delta_1, \gamma_1 \text{ and } \gamma_2$, for which the joint-distribution is to be estimated from data (see SI for model details). We use $\alpha = 1.1, \beta = 0.08, \delta_1 = 0.01, \delta_2 = 0.002, \delta_3 =0.002, \gamma_1=0.1, \gamma_2 = 0.1 \text{ and } \gamma_3 =0.06 $ as the true parameter values. Using the data for only first 15 days, we follow the same sampling procedure as described previously for the SEIR model. 
The marginal distributions along with their mean and mode are shown in fig. \ref{fig:mcmc_dist_sepihr}. 

\begin{figure}[h!]
\vspace{-0.25cm}
    \centering
    \includegraphics[width = \linewidth]{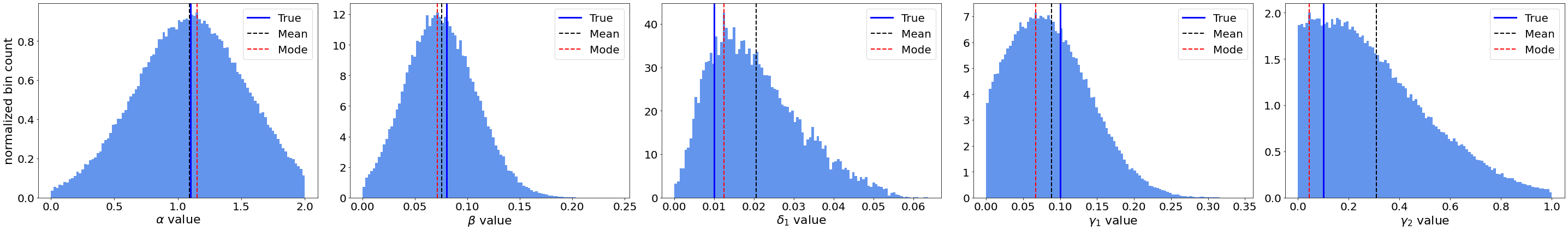}
    \vspace{-0.45cm}
    \caption{SEPIHR: Empirical distribution after $3\times 10^5$ samples from the MCMC sampling procedure.}
  \label{fig:mcmc_dist_sepihr}    
  \vspace{-0.25cm}
\end{figure}
 We note that mode is very close to the true values in both the models,  thus validating the capability of the sampling procedure in correctly estimating the parameter values.
\vspace{-0.2cm}
\paragraph{Related literature:} 
Before concluding this section, we briefly mention related literature.
Variational inference (VI) based approach of \cite{nico} provides improvements over \cite{dondelinger}, however due to modelling assumptions is not suited for our work. The optimization based gradient matching approaches of \cite{ramsay,hua,gonzalez,niu,odin} and others like \cite{root_n_estimator}  only provide point-estimates. The generative modelling approach of  \cite{barber} suffers from identifiablity issues as explained by \cite{controversy}.
Approaches with different sampling methods would include  \cite{hmc_efficient, emulation_acc_hmc, aistats_rev2_paper2, thermodynamic_integration,geometric_sampling} and approximation based methods would include \cite{abc_param_inference,aistats_rev2_paper1,laplace_based_ode_parameter}.
Other VI based methods would include \cite{meeds, ghosh21b}.
Probabilistic numerics \cite{probablistic_numerics} based methods include \cite{ fenrir} and \cite{fast_inversion_likelihood}.
\cite{teymur_implicit_prob_integrators} showed the use of  probabilistic integrators for ODEs in parameter estimation.
\cite{bayesian_soln_UQ_DE,aistats_rev2_paper3,aistats_rev2_paper4,sampling_stan} are other useful references.
\vspace{-0.25cm}
\section{Tractable Scenario-set Construction} \label{sec:sce_reduction}
\vspace{-0.35cm}
In section \ref{sec:empirical_sampling_result}, we obtained an empirical joint distribution on the SEIR parameters $\alpha, \beta, \gamma$ in the form of $3\times 10^5$ samples. Each of these samples represent a real-world scenario and will be used to represent uncertainty through the scenario-set in the vaccine allocation stochastic optimization formulation (discussed in the next section). However, working with such large number of samples is  computationally prohibitive, therefore we first discuss how to reduce the number of these samples (or scenarios) while still correctly representing the joint-distribution. We introduce some mathematical preliminaries:

Let $\mathbb{P} = \sum_{i \in I} p_i\delta_{\xi_i}$, where $\xi_i \in \mathbb{R}^d$ represents the location and $p_i \in [0,1]$ represents the probability of the $i$-th scenario in   $\mathbb{P}$, where $i \in I = \{1,\dots,n\}$. 
Let $\mathbb{Q}$ represents the target distribution. 
$\mathbb{Q} = \sum_{j \in J} q_j\delta_{\zeta_j}$, where $\zeta_j \in \mathbb{R}^d$ represents the location and $q_j \in [0,1]$ represents the probability of the $j$-th scenario in   $\mathbb{Q}$, where $j \in J =\{1,\dots,m\}$. 

We can now define the type-l Wasserstein distance between original $\mathbb{P}$ and target $\mathbb{Q}$ as following:
\begin{alignat}{10}
d_{l}(\mathbb{P}, \mathbb{Q}) &= \min_{\pi \in \mathbb{R}_{+}^{n \times m }} \Big( \sum_{i \in I}  \sum_{j \in J} \pi_{ij} ||\xi_i - \zeta_j  ||^l \Big)^{1/l} \label{eq:wass_dist} \\
& \phantom{} \text{ s.t. } \sum_{j \in J} \pi_{ij} = p_i \forall i \in I \;;\; 
\sum_{i \in I} \pi_{ij} = q_j \forall j \in J  \notag 
\end{alignat}
where $l \ge 1$. The linear program \eqref{eq:wass_dist} corresponds to the min-cost transportation problem, where $\pi_{ij}$ represents the amount of mass moved from location $\xi_i$ to $\zeta_j$ and $||\xi_i - \zeta_j||^l$ represents the associated cost incurred in moving unit mass.

Let $\mathcal{P}_{u}(\mathcal{R},n)$ denote the set of all \emph{uniform} discrete distributions  with exactly $n$ scenarios  on any given space $\mathcal{R}$ where $ \mathcal{R} \subseteq \mathbb{R}^d$. Similarly, let 
$\mathcal{P}(\mathcal{R},m)$ denote the set of all discrete distributions (not necessarily uniform)  with \emph{at-most} $m$ scenarios.

Let the  discrete probability distribution over parameters (for eg: SEIR $\alpha, \beta, \gamma$) obtained using sampling be denoted as $\hat{\mathbb{P}}$, where $\hat{\mathbb{P}}$ belongs to the set  $\mathcal{P}_{u}(\mathcal{R},n)$, i.e. $\hat{\mathbb{P}} \in \mathcal{P}_{u}(\mathcal{R},n)$, where $\mathcal{R}$ is determined by the upper and lower bounds on $\alpha, \beta, \gamma$ and $n = 3\times 10^5$.\\
The scenario reduction problem can now be defined as:
\begin{equation}\label{eq:sc_reduction}
\mathcal{CSR}(\hat{\mathbb{P}},m  ) = \min_{\mathbb{Q} \in \mathcal{P}(\mathcal{R}, m) } d_{l}(\hat{\mathbb{P}}, \mathbb{Q})
\end{equation}

\begin{figure}[h!]
\vspace{-0.25cm}
    \centering
    \includegraphics[width = \linewidth]{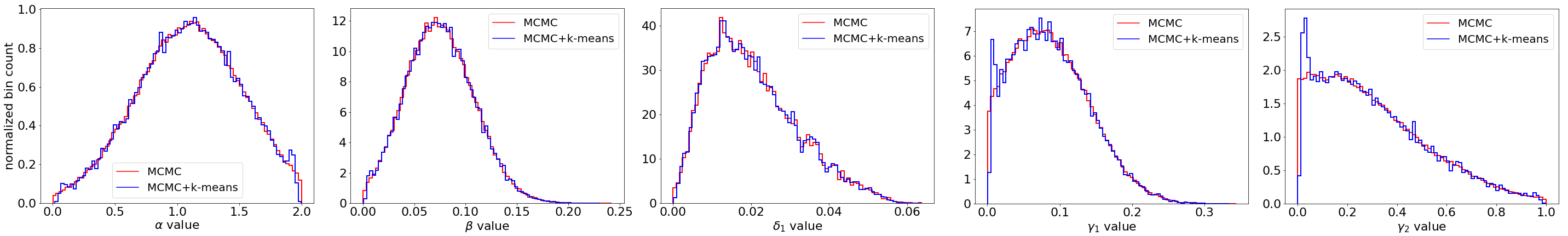}
    \caption{SEPIHR: MCMC+k-means denotes the empirical distribution with $3\times 10^4$ samples obtained after doing k-means clustering on the original $3\times 10^5$ MCMC samples shown in fig. \ref{fig:mcmc_dist_sepihr}.}
    \label{fig:kmeans_sepihr}
    \vspace{-0.3cm}
\end{figure}

Solving $\mathcal{CSR}$ defined in \eqref{eq:sc_reduction} is not easy for an arbitrary $l$. Fortunately for our purposes, the $\mathcal{CSR}$ problem is same as the k-means clustering problem with $m=\text{k}$ and $l=2$ \cite{rujeerapaiboon_sc_reduction}. 
Clustering distributes the set of $n$ points in $I$ with locations $\xi_i $  into \text{k} mutually exclusive subsets $I_1, \dots I_k$.  The $\zeta_j $'s (also known as centroids) and their associated probability $q_j$'s of the target distribution $\mathbb{Q}$ are given by:
\begin{wrapfigure}{r}{0.6\linewidth}
\vspace{-0.25cm}
  \centering
      \captionsetup{singlelinecheck=on, margin={0.15cm, 0.45cm}, justification=justified, format=plain}
  \includegraphics[width=\linewidth]{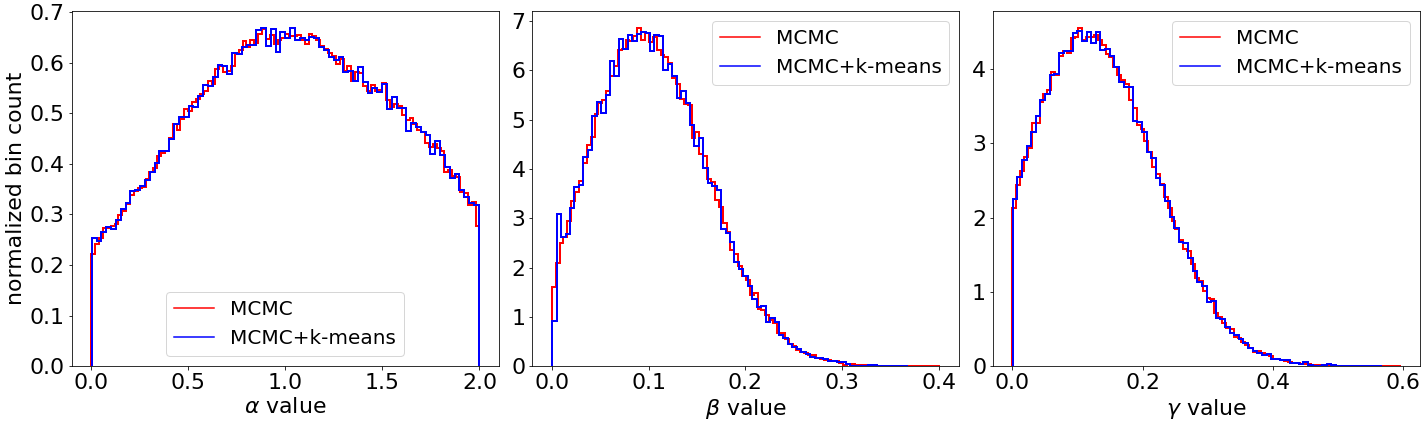}
  \vspace{-0.4cm}
  \captionof{figure}{SEIR: MCMC+k-means denotes the empirical distribution with $3\times 10^4$ samples obtained after doing k-means clustering on the original $3\times 10^5$ MCMC samples.}
  \label{fig:kmeans}
  \vspace{-0.2cm}
\end{wrapfigure}
\begin{equation}\label{eq:target_Q}
\begin{small}
\hspace{-0.5cm}\displaystyle \zeta_j = \frac{1}{|I_j|}\sum_{i \in I_j}{\xi_i} \text{ and } q_j = \frac{|I_j|}{n}.
\end{small}
\end{equation}
 Although theoretically, \text{k}-means clustering is known to be NP-hard \cite{mahajan_kmeans_nphard, aloise_kmeans_nphard}, however empirically a high quality solution can be easily obtained using Lloyd's algorithm \cite{lloyd_kmeans}  also commonly known as  the  k-means clustering algorithm.
We run k-means clustering on the original $3\times10^5$ samples with $\text{k}=3\times 10^4$ to get a 10x reduction. For the SEIR model, resulting distribution (denoted $\hat{\mathbb{Q}}$), marginal distributions are shown in fig. \ref{fig:kmeans}. We observe that the reduced target distribution is very close to the original distribution. Using  $\hat{\mathbb{Q}}$, we can now construct our scenario-set (denoted $\Omega$) as:
\begin{equation} \label{eq:def_Omega}
\Omega  = \big \{ (\alpha^j, \beta^j, \gamma^j ,p_j) \; \forall \;j \in  \{1,\dots,\text{k} \}   \big \}
\end{equation}
where $(\alpha^j, \beta^j, \gamma^j)  = \hat{\zeta}_j \text{ and } p_j= \hat{q}_j$. 
As $\text{k} = 3\times 10^4$, therefore $\Omega$  also has $3\times 10^4$ 
scenarios, i.e. $|\Omega| = 3\times 10^4$. We will henceforth use this $\Omega$ as our scenario set. The mode of $\hat{\mathbb{Q}}$ corresponds to the nominal estimate of $\alpha,\beta \text{ and } \gamma$. Similarly, the result for SEPHIR model (fig. \ref{fig:mcmc_dist_sepihr}) is shown in fig. \ref{fig:kmeans_sepihr}. Experiments with values of k, other than $\text{k}=3\times 10^4$, for both models are provided in SI.

\vspace{-0.35cm}
\section{Optimal Vaccine Allocation Formulation and Solution Algorithm} \label{sec:vacc_allocation}
\vspace{-0.3cm}
We now work towards formulating our optimization problem for vaccine allocation. Our goal is to allocate vaccines (on a daily basis) to a set of $\mathbb{K}$ sub-populations, such that the maximum number of total infections is minimized. This objective ensures that the peak of the pandemic is minimised as much as possible in order to reduce the burden on the healthcare services particularly medical personnel at the height of the pandemic. The $\mathbb{K}$ sub-populations correspond to different geographical regions such as nearby cities in a state. Let $\mathcal{K} = \{1,\dots,\mathbb{K} \}$.

\begin{figure}[h!]
\vspace{-0.6cm}
\begin{minipage}{0.24\linewidth}
   \centering 
\includegraphics[scale = 0.35]{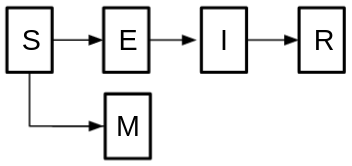}
 \caption{SEIRM model}\label{fig:seirm}
\end{minipage}\hfill
\begin{minipage}{0.73\linewidth}
\vspace{-0.35cm}
\begin{adjustwidth}{-0.4cm}{-0.2cm}
\begin{small}
\begin{equation}\label{eq:seirm}
\begin{rcases}
\begin{alignedat}{10}
&\frac{dS_k(t)}{dt}   &&\coloneqq \dot{S}_k(t) = - \eta V_k(t)   - \frac{u_k(t) \alpha}{N_k} \big(S_k(t)- \eta V_k(t)\big)\big(\sum_{r=1}^{\mathbb{K}} {\lambda}_{r}^{k}I_r(t)  \big)  \\[-1em]
&\frac{dE_k(t)}{dt}   &&\coloneqq \dot{E}_k(t) =  \frac{u_k(t)\alpha}{N_k} \big(S_k(t)- \eta V_k(t)\big)\big(\sum_{r=1}^{\mathbb{K}} {\lambda}_{r}^{k}I_r(t)  \big) -  \beta E_k(t)   \\[-0.8em]
&\frac{dI_k(t)}{dt}   &&\coloneqq \dot{I}_k(t) = \beta E_k(t) - \gamma  I_k(t)     \\
&\frac{dR_k(t)}{dt}   &&\coloneqq \dot{R}_k(t) = \gamma I_k(t), \phantom{MMmmm}\;
\frac{dM_k(t)}{dt}   \coloneqq \dot{M}_k(t) = \eta V_k(t) 
\end{alignedat}
\end{rcases}\hspace{-0.25cm}
\end{equation}
\end{small}
\end{adjustwidth}
\end{minipage}
\vspace{-0.35cm}
\end{figure}

The spread of disease in each sub-population is modeled using a separate SEIR model.
To account for the vaccinated individuals, the SEIR model in fig. \ref{fig:seir} is updated with a new compartment (denoted by M) to represent the immune population and the updated model (fig. \ref{fig:seirm}) is denoted by SEIRM. Let $V_k(t)$ represent the number of people vaccinated at time $t$ in the $k$-th sub-population  and $\eta$ be efficacy of the vaccine, then the ODEs corresponding to the SEIRM model of the $k$-th sub-population are given by eq: \eqref{eq:seirm}.

We also consider two important features of disease transmission. First, due to mobility there is contact between infected individuals of one sub-population with the susceptible individuals of another sub-population. Second, due to different levels of mobility between different sub-populations,
onset of the pandemic in each of the sub-populations generally vary. Both of these are accounted in the updated states $S_k$ and $E_k$ in eq:\eqref{eq:seirm},
where 
${\lambda}_{r}^{k}$ denotes the mobility levels from sub-population $r$ to sub-population $k$ and $u_k(t)$ corresponds to a sigmoid function, $u_k(t) \coloneqq 1/(1 + e^{-c_1^k(t-c_2^k)} )$ with parameters $c_1^k \text{ and } c_2^k$.
In particular, $c_2^k$ controls the onset of the pandemic in the $k$-th sub-population, therefore  we also account for uncertainty in $c_2^k \; \forall\; k \in \mathcal{K}$, by appropriately extending the scenario set $\Omega$ (for details see SI).

Let $\mathcal{T} = \{1,\dots,T\}$, denote the simulation time period, $\mathcal{T}_v = \{t_s,\dots,t_l\}$ denote the vaccination time-period   where  $t_s \text{ and } t_l$ are the first and last  vaccination days, such that $\mathcal{T}_v \subseteq \mathcal{T}$. We can now write the \emph{nominal (or non-stochastic)} optimization problem (denoted $\mathcal{NF}$) for vaccine allocation as \eqref{eq:opt_prob_nom},
where $B_t$ in \eqref{eq:nominal_total_B} denotes the total daily vaccine budget for all $\mathbb{K}$ sub-populations and $U_t^k$ in \eqref{eq:nominal_ind_B} denotes the vaccine budget for $k$-th sub-population. Equations \eqref{eq:nominal_odes} represent the ODE constraints, \eqref{eq:nominal_Isum} \& \eqref{eq:nominal_Imax} together computes the maximum (or peak) infection of the total population (denoted $\mathcal{I}$) and \eqref{eq:nominal_obj} minimizes the peak infection.

We now provide the \textbf{uncertainty-informed}, i.e. stochastic counterpart (denoted $\mathcal{SF}$) of the nominal problem $\mathcal{NF}$ in \eqref{eq:opt_prob_stoc}, where $\Omega$ denotes the scenario-set, recall \eqref{eq:def_Omega}. Each state S,E,I,R,M in ODE constraints in \eqref{eq:stoc_sc}   now has an associated superscript $\omega$ corresponding to that scenario,  
$\mathcal{I}_{\omega}$ denotes the peak infection for scenario $\omega$, \eqref{eq:stoc_obj} computes the expected peak infection over all scenarios. The vaccine budget constrains in \eqref{eq:stoc_budget_v} remain same as in $\mathcal{NF}$.

\begin{adjustwidth}{-0.25cm}{-0.05cm}
\begin{minipage}[t]{0.42\linewidth}
\vspace{-0.3cm}
\begin{subequations}\label{eq:opt_prob_nom}
\begin{alignat}{10}
&  \mathcal{NF}: \min_{V} \;\; \mathcal{I} \qquad \qquad \text{(Nominal)} \label{eq:nominal_obj} \\
& \text{s.t.} \notag \\
&  \{\text{\eqref{eq:seirm}}\} \phantom{MMMM}\; \forall \;\;  k \in \mathcal{K}, t \in \mathcal{T}  \label{eq:nominal_odes} \\
& \small  \sum_{k=1}^{\mathbb{K}} I_k(t)  = \tilde{I}(t) \phantom{Mmmm} \forall \;\; t \in \mathcal{T} \label{eq:nominal_Isum}\\
&  \small \tilde{I}(t)  \le \mathcal{I} \phantom{MMMMmmm}\; \forall  \;\; t \in \mathcal{T} \label{eq:nominal_Imax}\\
&   \sum_{k=1}^{\mathbb{K}} V_k(t)  \le B_{t} \phantom{mmmm}\; \forall \;\;  t \in \mathcal{T}_v \label{eq:nominal_total_B}\\
& 0 \le V_{k}(t) \le U_{t}^k  \phantom{m} \forall \; k \in \mathcal{K}, t \in \mathcal{T}_v \label{eq:nominal_ind_B} \\
&  V_{k}(t) = 0  \quad\; \forall \;\; k \in \mathcal{K}, t \in \mathcal{T} \setminus  \mathcal{T}_v \label{eq:nominal_v_0}
\end{alignat}
\end{subequations}
\end{minipage} \hspace{0.2cm}\vline\hspace{0.55cm}
\begin{minipage}[t]{0.52\linewidth}
\vspace{-0.3cm}
\begin{adjustwidth}{-0.45cm}{0cm}
\begin{subequations}\label{eq:opt_prob_stoc}
\begingroup
\begin{alignat}{10}
&\hypertarget{SF}{\mathcal{SF}}: \min_{V} \quad \sum_{\omega \in \Omega } p_{\omega}\mathcal{I}_{\omega} \phantom{MMMM}\; \text{(Stochastic)}  \label{eq:stoc_obj} \\
& \color{blue} \begin{rcases} 
\scriptstyle
\color{black} \dot{S}_k^{\omega}(t)   = - \eta V_k(t)   - \frac{u^{\omega}(t)\alpha^{\omega}}{N} \big(S_k^{\omega}(t)- \eta V_k(t)\big)\big(\sum_{r=1}^{\mathbb{K}} \lambda_{r}^k I_r^{\omega}(t)  \big)      \\
\scriptstyle \color{black} \dot{E}_k^{\omega}(t)   =  \frac{u^{\omega}(t)\alpha^{\omega}}{N} \big(S_k^{\omega}(t) - \eta V_k(t)\big)\big(\sum_{r=1}^{\mathbb{K}} \lambda_{r}^kI_r^{\omega}(t)  \big) -  \beta^{\omega} E_k^{\omega}(t)  \\
 \end{rcases}   \notag\\
&{\color{blue} \begin{drcases}  \color{black}
\dot{I}_k^{\omega}(t)   =  \beta^{\omega} E_k^{\omega}(t) - \gamma^{\omega}  I_k^{\omega}(t)     \\
\color{black} \dot{R}_k^{\omega}(t)   =  \gamma^{\omega} I_k^{\omega}(t) ,\;
\color{black} \dot{M}_k^{\omega}(t)    = \eta V_k(t) \\
\color{black} \sum_{k=1}^{\mathbb{K}} I_k^{\omega}(t)  = \tilde{I}^{\omega}(t) \phantom{M} , \color{black} \tilde{I}^{\omega}(t)  \le \mathcal{I}_{\omega}  \\ 
  \end{drcases} \bigforall \; 
  \begin{matrix*}[l]
  k \in \mathcal{K},  \\
  t \in \mathcal{T},  \\ 
  \omega \in \Omega
  \end{matrix*}
    } \hspace{-0.2cm}   \label{eq:stoc_sc}\\
 & \{\text{\eqref{eq:nominal_total_B}, \eqref{eq:nominal_ind_B}, \eqref{eq:nominal_v_0}} \} \label{eq:stoc_budget_v} 
\end{alignat}
\endgroup
\end{subequations}
\end{adjustwidth}
\end{minipage}
\end{adjustwidth}

\begin{definition} A vaccine policy $\mathcal{V}$ is defined as: $\mathcal{V} =\{ V_k(t) \; \forall \;  k \in \mathcal{K}, t \in \mathcal{T} \}.$
\end{definition}

\begin{theorem}\label{th:feasibiltiy}{Feasibility of  $\; \mathcal{V}$:}
The feasibility of a vaccine policy $\mathcal{V}$ in $\mathcal{SF}$ is only decided by the budget constraints in \eqref{eq:stoc_budget_v} and not by the ODE constraints in \eqref{eq:stoc_sc}.
\end{theorem}
Theorem \ref{th:feasibiltiy} holds true because of the budget constraints \eqref{eq:stoc_budget_v}, $V_k(t)$ is non-negative and finite. Thus the existence and uniqueness of a solution to ODEs in \eqref{eq:stoc_sc}  is guaranteed and can be shown analytically using the Picard-Lindel{\"o}f theorem with appropriate initial conditions \cite{levinson, nonlinear-systems-sastry,existence_uniqueness1}.

\begin{lemma}\label{lemma:decomposability}{Decomposability w.r.t $\Omega$:}
 For a given (fixed) vaccine policy $\mathcal{V}$, the ODE constraints in \eqref{eq:stoc_sc} become decomposable, i.e. the set of ODE constraints in scenario $\omega_i $ can be solved independently of the set of ODE constraints in scenario $\omega_j \; \forall j \in \Omega \setminus {i}$.
 \end{lemma}
Lemma \eqref{lemma:decomposability} follows from the fact that for a given scenario (say $\omega_i$), constraints in \eqref{eq:stoc_sc} require parameters only corresponding to scenario $\omega_i$. This has major computational implications as it allows for parallel evaluation of scenarios in $\Omega$. Due to the additive nature of the objective function  \eqref{eq:stoc_obj} w.r.t. to $\Omega$, we can compute the objective function value after parallel computation of scenarios. 
Therefore, we can efficiently solve $\mathcal{SF}$ using an iterative heuristic based optimization procedure described in algorithm \ref{algo:1}. Details on heuristics are provided in SI.

\vspace{-0.1cm}
\begin{algorithm}[H]
\caption{Optimization procedure to solve $\mathcal{NF}$ or $\mathcal{SF}$}\label{algo:1}
  \footnotesize
\begin{algorithmic}[1]
\STATE Randomly sample a batch of Vaccine policies of size $B$, i.e. $\bar{\mathcal{V}}_0 = \{ \mathcal{V}_1, \dots, \mathcal{V}_B \}$ and set $i = 0$.

\WHILE{$i \le N_{\text{opt}} $}
\FOR{$k \gets 1$ to $B$}                    
    \STATE Evaluate constraint violation (denoted $C_k$) 
    of $\bar{\mathcal{V}}_i[k]$ using \eqref{eq:stoc_budget_v} .
    
    \STATE In \emph{parallel}, simulate all $|\Omega|$ scenarios for $\bar{\mathcal{V}}_i[k]$ using an ODE solver to compute $\mathcal{I}_{\omega} $.  
    
    \STATE Compute $f_{obj}^k$\eqref{eq:stoc_obj}:  $f_{obj}^k \gets \sum_{\omega \in \Omega} p_{\omega}  \mathcal{I}_{\omega} $ 
\ENDFOR
\STATE Update the batch of vaccine policies  $\bar{\mathcal{V}}_{i}$ with heuristic rules  using $\{f_{obj}^1,\dots,f_{obj}^B\}$ and $\{C_1, \dots,C_B \} $ to generate next batch of vaccine policies  $\bar{\mathcal{V}}_{i+1}$.
\STATE $i \gets i+1$ 
\ENDWHILE \\
\RETURN feasible vaccine policy $\mathcal{V}$ with lowest $f_{obj}$.
\end{algorithmic}
\end{algorithm}

\vspace{-0.5cm}
\section{Experimental (Simulation) Results}\label{sec:exp_simu_results}
\vspace{-0.3cm}
We show the efficacy of our proposed approach on two different disease transmission models, i.e. the SEIR and the SEPIHR models. 
For all experiments we report average of 5 runs. In addition to the experiments in this section, various other numerical experiments under different setups are provided in SI.

Recall that in section \ref{sec:empirical_sampling_result},  we have already discussed the details and sampling results for both the SEIR and SEPIHR model including figures \ref{fig:mcmc_dist} \& \ref{fig:mcmc_dist_sepihr} respectively. Also, in section \ref{sec:sce_reduction} (including fig. \ref{fig:kmeans} \& \ref{fig:kmeans_sepihr}), we have discussed how to obtain a tractable scenario set $\Omega$ using k-means to account for uncertainty in the vaccine allocation. Therefore our main goal in this section is to show the benefit of incorporating uncertainty by comparing the vaccine allocation policy (denoted  $\mathcal{V}_{\mathcal{N}}$) obtained from solving the nominal formulation $\mathcal{NF}$   against the vaccine allocation policy (denoted  $\mathcal{V}_{\mathcal{S}}$) obtained from solving the stochastic solution $\mathcal{SF}$ . We also benchmark against a zero or no-vaccination policy denoted (denoted  $\mathcal{V}_{\phi}$), where $\mathcal{V}_\phi = \{V_k(t) =0 \;\forall\; t \in \mathcal{T}, k \in \mathcal{K} \}$.

\textbf{SEIR model}:  We use a total simulation time horizon of $T=120$ days, vaccination period of 25 days starting on $t_s =16$ and ending on $t_l=40$ with daily available vaccine budgets $B_t = 24 \times 10^3$ and $U_t^k = 10^4$. Importantly, note that in section \ref{sec:empirical_sampling_result} for parameter estimation, we used data only for first 15 days, i.e. $T=\{1,\dots,15\}$, thus maintaining consistency for real-world applicability. 
We perform experiments in two different settings, in the \textbf{first setting} we work with $\mathbb{K}=3$ i.e. three sub-populations of sizes $7.5\times 10^5, 5\times 10^5 \text{ and } 10^6$ respectively and in the \textbf{second setting} we increase $\mathbb{K}$ to $\mathbb{K}=4$, with an additional sub-population of size $6\times 10^5$. Numerical values of other parameters like $\lambda_r^k,c_1^k,c_2^k, \eta$ and additional experiments to evaluate their effect are provided in SI. 

For each setting i.e. $\mathbb{K}=3$ and $\mathbb{K}=4$, using algorithm \ref{algo:1} and the nominal estimates of $\alpha,\beta \text{ and } \gamma$, we solve the $\mathcal{NF}$ to get the nominal vaccine policy $\mathcal{V}_{\mathcal{N}}$. Using the scenario-set $\Omega$ (generated from the discrete-parameter distribution) we solve $\mathcal{SF}$ to get the uncertainty-informed vaccine policy $\mathcal{V}_{\mathcal{S}}$.\\ 
We next evaluate the efficacy of all the three vaccine policies i.e. $\mathcal{V}_\phi, \mathcal{V}_{\mathcal{N}} \text{ and } \mathcal{V}_{\mathcal{S}}$. For each of these policies,  we simulate all the scenarios in the scenario-set $\Omega$  and compute the expected values of all the states i.e. S,E,I,R and M over the time horizon $T$.

\begin{figure}[h!]
\begin{adjustwidth}{-0.4cm}{-0.4cm}
\centering
\captionsetup[subfigure]{aboveskip=-1pt,belowskip=-1pt}
\begin{subfigure}{.5\linewidth}
  \centering
  \includegraphics[width=\linewidth]{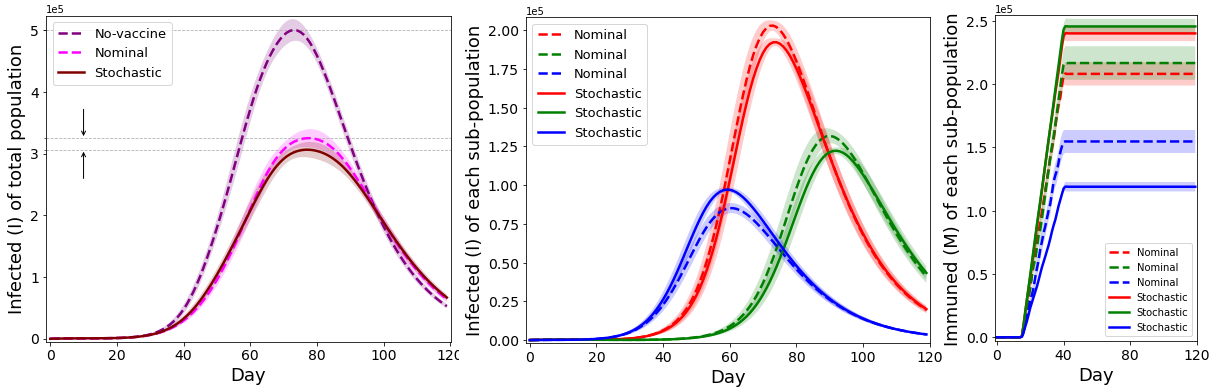}
  \caption{$\mathbb{K}=3$ sub-populations}
  \label{fig:seirm_pop3}
\end{subfigure}%
\begin{subfigure}{.5\linewidth}
  \centering
  \includegraphics[width=\linewidth]{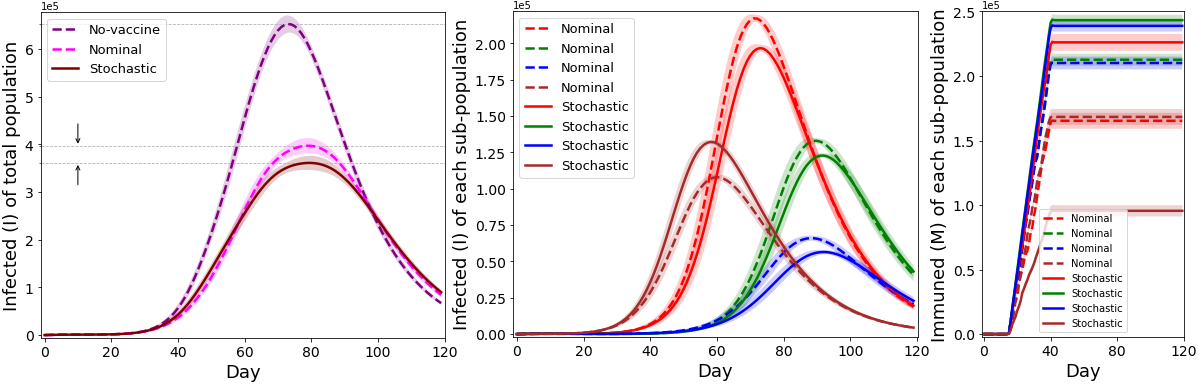}
  \caption{$\mathbb{K}=4$ sub-populations}
  \label{fig:seirm_pop4}
\end{subfigure}
\vspace{-0.5cm}
\caption{\textbf{SEIR:} Evaluation of different vaccine policies i.e. no-vaccine $\mathcal{V}_{\phi}$, nominal $\mathcal{V}_{\mathcal{N}}$ and stochastic $\mathcal{V}_{\mathcal{S}}$.}
\label{fig:seirm_exps}
\end{adjustwidth}
\vspace{-0.4cm}
\end{figure}

The evolution of the infected state (I)  of the \emph{total population} and the infected (I) and immuned (M) states of each \emph{sub-population} are shown in fig. \ref{fig:seirm_pop3}  and \ref{fig:seirm_pop4} for $\mathbb{K}=3$ and $4$ respectively. We note that for $\mathbb{K}=3$ (fig \ref{fig:seirm_pop3}) ,  the expected peak infection is reduced from around  501k (with no-vaccination i.e. $\mathcal{V}_{\phi}$) to 324k with nominal vaccination policy $\mathcal{V}_{\mathcal{N}}$. This reduction of peak infection by 35.3\% is expected due to vaccination. More importantly, we observe that with the stochastic vaccination policy $\mathcal{V}_{\mathcal{S}}$ the peak infection is further reduced to around 308k, which is an improvement of around \textbf{4.9\%} over $\mathcal{V}_{\mathcal{N}}$ and 38.56\% over $\mathcal{V}_{\phi}$. This improvement of $\mathcal{V}_{\mathcal{S}}$ over $\mathcal{V}_{\mathcal{N}}$ by \textbf{4.9\%} is also referred to as the \emph{value of stochastic solution} (VSS) or equivalently the benefit of accounting for uncertainty.

For $\mathbb{K}=4$ (fig \ref{fig:seirm_pop4}), we observe that the peak infection with no-vaccine policy $\mathcal{V}_{\phi}$ is around 653k, and is reduced to 393k with $\mathcal{V}_{\mathcal{N}}$ and is further reduced to 361k with $\mathcal{V}_{\mathcal{S}}$, i.e. $\mathcal{V}_{\mathcal{S}}$ provides a reduction of around \textbf{8\%} over $\mathcal{V}_{\mathcal{N}}$. This higher VSS of \textbf{8\%} for $\mathbb{K}=4$ compared to \textbf{4.9\%} for $\mathbb{K}=3$ is due to the fact that
the size of scenario set $|\Omega|$  is directly proportional to the number of sub-populations $\mathbb{K}$. Recall that we also account for the uncertainty in the onset of the pandemic in each sub-population through the parameter  $c_2^k$.

Note that since the immuned  sub-population size  is directly proportional to vaccines allocated to that sub-population, therefore the third figure in \ref{fig:seirm_pop3} and \ref{fig:seirm_pop4} also shows how many vaccines are allocated to each sub-population relative to each other. We observe that there is a clear difference between the nominal and the stochastic allocations. This significant difference in nature of the vaccine policies explain the reduction of \textbf{4.9\%} and \textbf{8\%} respectively, providing validity to our results in the sense that the reductions obtained are not simply due to minor numerical changes in solution values. We further discuss the differences of the two policies ($\mathcal{V}_{\mathcal{N}}$ vs $\mathcal{V}_{\mathcal{S}}$)  in SI. 

\textbf{SEPIHR model:} We next evaluate our approach on the SEPIHR model with additional states P (for protective quarantine) and H (for hospitalised quarantined).
Importantly here as the number of hospitalisations (H) is  modeled explicitly, therefore we minimize the peak (maximum) hospitalisations. For this model, corresponding vaccine allocation optimization formulations (i.e. nominal and stochastic) are provided in SI.

\begin{figure}[h!]
\vspace{-0.3cm}
\centering
\begin{adjustwidth}{-0.3cm}{-0cm}
\begin{subfigure}{0.225\linewidth}
  \centering
  \includegraphics[width=0.98\linewidth]{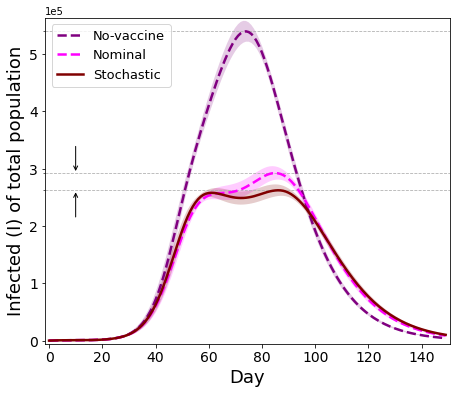}
  \label{fig:rank_1_tot_I}
\end{subfigure}%
\begin{subfigure}{0.225\linewidth}
  \centering
  \includegraphics[width=\linewidth]{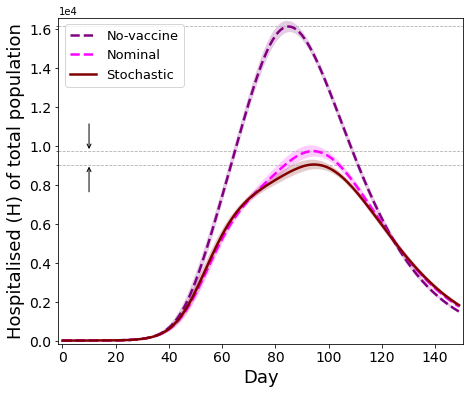}
  \label{fig:rank_1_tot_H}
\end{subfigure}%
\begin{subfigure}{0.225\linewidth}
  \centering
  \includegraphics[width=\linewidth]{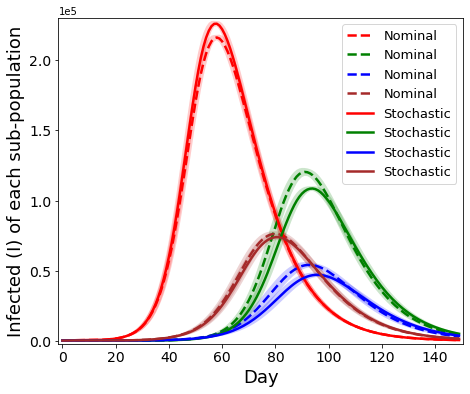}
  \label{fig:rank_1_ind_I}
\end{subfigure}%
\begin{subfigure}{0.225\linewidth}
  \centering
  \raisebox{0.45cm}{\includegraphics[width=0.96\linewidth]{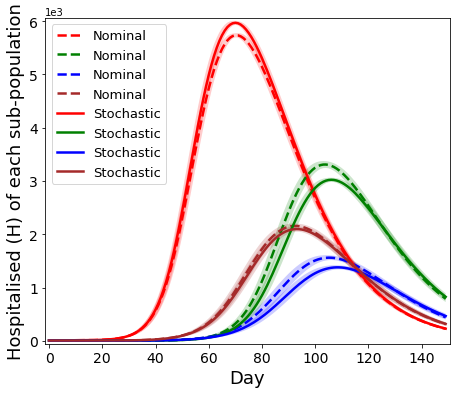}}
  \label{fig:rank_1_ind_H}
 \end{subfigure}%
 \begin{subfigure}{0.13\linewidth}
   \centering
  \includegraphics[width = \linewidth]{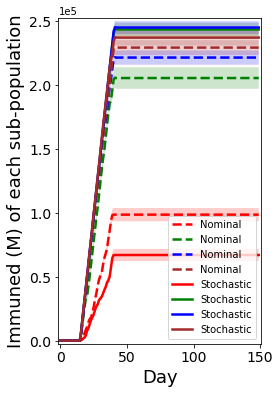}
   \label{fig:rank_1_ind_M}
 \end{subfigure}
 \vspace{-1cm}
\caption{\textbf{SEPIHR:} Evaluation of different policies: $\mathcal{V}_{\phi}$, $\mathcal{V}_{\mathcal{N}}$ and $\mathcal{V}_{\mathcal{S}}$ with $\mathbb{K}=4$ sub-populations.} \label{fig:sepihrm_exps}
\end{adjustwidth}
\vspace{-0.3cm}
\end{figure}

In fig. \ref{fig:sepihrm_exps}, we show the evolution of the infected (I) and hospitalised (H) states of the total population, along with the  I,H and immuned (M) for each of the 4 sub-populations. We note that the peak infections (I) for the three policies (i.e. $\mathcal{V}_{\phi}$, $\mathcal{V}_{\mathcal{N}}$ and $\mathcal{V}_{\mathcal{S}}$)  are around  539k,	280k  and 262k  and the peak hospitalisations are around 16k, 9.4k and 9k respectively. Therefore, $\mathcal{V}_{\mathcal{S}} $ provides a reduction of \textbf{6.3\%} in peak infections (I) over $  \mathcal{V}_{\mathcal{N}}$ and a \textbf{4.4\%} reduction in peak hospitalisations (H). Interestingly, from the fifth plot in fig. \ref{fig:seirm_exps}, we note that despite its largest size and earliest the onset of the pandemic, red population is allocated the least vaccines. This can be explained by the fact that we aim to minimize the peak of the total population. In SI, we provide more such sub-population level discussions on the differences in nature the of the optimal policies including $\mathcal{V}_{\mathcal{N}}$ vs $\mathcal{V}_{\mathcal{S}}$.

The  above results on SEIR and SEPIHR models clearly demonstrate the benefit of uncertainty-informed vaccine allocation using Bayesian inference over using nominal estimates.
Our improvements of \textbf{4-8\%}  are either consistent with prior works in literature such as \cite{vacc_2phase_stochastic} or much better \cite{thul_powell}.\\
In SI, we provide more experiments under different setups for both the models. We also discuss the possible societal impact of our work.

\vspace{-0.3cm}
\section{ Concluding Remarks and Future Work}\label{sec:concluding_remarks}
\vspace{-0.35cm}
In this paper, we proposed an uncertainty informed vaccine allocation problem as a stochastic optimization problem, for which the tractable scenario-set is constructed in a novel data-driven manner using Bayesian inference for ODEs with GPs. We also proposed a scalable solution algorithm to solve the stochastic program and showed that a significant gain can be achieved by accounting for uncertainty. For future work, a natural extension would be to systematically investigate equity and fairness of allocation through additional constraints and different objective functions.

\setcounter{section}{0}

\bibliographystyle{plain}
\bibliography{paper}

\newgeometry{left=1in, right=1in,
    textheight=9in,
    top=1in,
    headheight=12pt,
    headsep=25pt,
    footskip=30pt}

\definecolor{egyptianblue}{rgb}{0.06, 0.2, 0.65}

\section*{\centering \Large Supplementary Information (SI)}

\section{Societal Impact}
Previous studies \cite{bertsimas_vaccine,acemoglu_parise,candogan,ghobadi}  have shown that accurate and data-informed  modeling of the spread of a pandemic is extremely useful for policy makers and public health agencies to make time-sensitive critical decisions. These decisions generally include imposing lockdowns, allocation of resources such as medical personnel, medical supplies (PPE kits, masks, oxygen etc.), testing facilities and \emph{vaccines}. 
Needless to mention that public health and economic recovery of our society crucially depends on these decisions.

Our approach and results suggest that accounting for uncertainty in key epidemiological parameters can improve the efficacy of time-critical allocation decisions.
We believe that our work will mainly have a positive societal impact.
However, at the same time it is extremely important to be mindful of the fairness and equity aspects of resource allocation. Given the flexible nature of our approach, we believe that our work can be leveraged to study and address these aspects systematically. 

\section{Definition of $\textbf{C}_{\phi_{k}}, \textbf{C}''_{\phi_{k}}, 
'\textbf{C}_{\phi_{k}} \text{ and } \textbf{C}'_{\phi_{k}} $ }

Let $\mathscr{K}(\cdot,\cdot)$ denote the a valid kernel function, then the entries of the matrices i.e. $\textbf{C}_{\phi_{k}}, \textbf{C}''_{\phi_{k}}, 
'\textbf{C}_{\phi_{k}} \text{ and } \textbf{C}'_{\phi_{k}} $ are given as:
\begin{alignat}{10}
\mathscr{K}(x_k(t), x_k(t')) &= \text{C}_{\phi_{k}}(t,t') \notag \\
\mathscr{K}(\dot{x}_k(t), x_k(t')) &= \frac{\partial \text{C}_{\phi_{k}}(t,t') }{\partial t}  \coloneqq\;&& \text{C}_{\phi_{k}}'(t,t') \notag \\
\mathscr{K}(x_k(t),\dot{x}_k(t')) &= \frac{\partial \text{C}_{\phi_{k}}(t,t') }{\partial t'}  \coloneqq\;&& '\text{C}_{\phi_{k}}(t,t') \notag \\
\mathscr{K}(\dot{x}_k(t),\dot{x}_k(t')) &= \frac{\partial^2 \text{C}_{\phi_{k}}(t,t') }{ 
\partial t  \partial t'}  \coloneqq\;\;&& \text{C}_{\phi_{k}}''(t,t') \notag 
\vspace{-0.2cm}
\end{alignat}

\section{GP hyper-parameter estimation}

For a zero-mean Gaussian process ($\bm{\mu}_k=\bm{0}$), with hyper-parameters $\bm{\phi}$ and $\bm{\sigma}$, the log-likelihood of $n$  observations (denoted $\textbf{y}$) at evaluation times $\textbf{t}$ can  be obtained in closed form as \cite{rasmussen}:

\begin{equation}\label{eq:likelihood}
\displaystyle \log(p(\textbf{y} | \textbf{t}, \bm{\phi}, \bm{\sigma}  )) = -\frac{1}{2}\textbf{y}^T( \textbf{C}_{\bm{\phi}}  + \bm{\sigma} \text{I} )^{-1}\textbf{y} -\frac{1}{2}\log|\textbf{C}_{\bm{\phi}}  + \bm{\sigma} \text{I} |
 -\frac{n}{2} \log 2 \pi
 \end{equation}

Equation \eqref{eq:likelihood} can be maximized w.r.t to $\bm{\phi}$
 and $\bm{\sigma}$. Note that the equation does not depend on the functional form of the ODEs.

\section{Accounting for uncertainty in $c_2^k$ }
The value of the parameter $c_2^k$ controls the onset of the pandemic in the $k$-th sub-population. We therefore also account for the uncertainty in the onset on the pandemic in each of the sub-populations.
We do this by augmenting the different scenarios in the scenario set $\Omega$. 
For each sub-population we consider $m$ different values of  $c_2^k$, where $m$ is around 4-5. Let $\Omega_{c_2^k} = \{  c_2^{k^1},\dots,c_2^{k^m}   \}$, then

$\Omega_{c_2} = \prod_{k=1}^{\mathbb{K}}\{  c_2^{k^1},\dots,c_2^{k^m}   \} $, or equivalently we can write:\\ $ \displaystyle \Omega_{c_2} = \big \{ \{ a_1,\dots, a_{\mathbb{K}} \} \;\; \forall \;\;  a_1  \in  \Omega_{c_2^1}, \dots ,  a_{\mathbb{K}}  \in  \Omega_{c_2^{\mathbb{K}}}  \big \} $ and $|\Omega_{c_2}| = |\Omega_{c_2^1}| \times \dots \times  |\Omega_{c_2^{\mathbb{K}}} | = m^{\mathbb{K}} $

Using $\Omega_{c_2}$,  we get the final scenario set as $\Omega_{\text{final}} = \Omega \times \Omega_{c_2}$, where $\Omega$ is the scenario set obtained after running k-means on the discrete parameter distribution. The probabilities of different scenarios in $\Omega_{\text{final}}$ can be  obtained easily by adjusting the probabilities of scenarios in $\Omega$ by a factor of $1/m^{\mathbb{K}}$.


\section{Numerical Values of different parameters}
For the experiments in the main text we used the following values:
We use $c_1^{k} =0.6 \; \forall \; k \; \in  \; \{1,2,3,4\}$ and $\eta = 0.99$. 
We use the following  base mobility matrices depending on the number of sub-populations:\\
\begin{minipage}{0.5\linewidth}
\begin{equation*}
\text{M}_3 = \begin{bmatrix}
    1        &  10^{-4}   &  10^{-4} \\
    10^{-4}  &  1         &  10^{-4} \\ 
    10^{-4}  &  10^{-4}   &  1    
\end{bmatrix} \text{ for } \mathbb{K}=3;
\end{equation*}    
\end{minipage}
\begin{minipage}{0.5\linewidth}
\begin{equation*}
\text{M}_4 = \begin{bmatrix}
    1        &  10^{-4}       &  10^{-4}    & 10^{-4} \\
    10^{-4}  &  1            &  10^{-4}    & 10^{-4} \\
    10^{-4}  &  10^{-4}      &  1          & 10^{-4} \\
    10^{-4}  &  10^{-4}      &  10^{-4}    & 1
\end{bmatrix} \text{ for } \mathbb{K}=4.
\end{equation*}  
\end{minipage}
$\lambda_r^k $ is equal to the $r,k$ entry of $\text{M}_{\mathbb{K}}$, i.e. $\lambda_r^k = \text{M}_{\mathbb{K}}[r,k]$.

\subsection{SEIR model}
For $\mathbb{K}=3$, we used 
$c_2^1 \in \{19, 20,21,22\}, c_2^2 \in\{29,30,31,32\},   c_2^3 \in \{9,10,11,12\} $ and the base mobility matrix $\text{M}_3$ with $\text{M}_3[0,2] = \text{M}_3[2,0] = 0$. \\
$B_t = 24 \times 10^3$ and $U_t^k = 10^4$.
For $\mathbb{K}=4$ we used 
 $c_2^1 \in  \{19, 20,21,22\}, c_2^2 \in\{29,30,31,32\}, c_2^3 \in \{24,25,26,27\},  c_2^4 \in\{9,10,11,12\} $ and the base mobility matrix $\text{M}_4$.\\
 $B_t = 32 \times 10^3$ and $U_t^k = 10^4$.
\subsection{SEPIHR model}
We used a total simulation time horizon of $T=150$ days, vaccination period of 25 days starting on $t_s =16$ and ending on $t_l=40$.
We used $c_2^1 \in\{9,10,11,12\}, c_2^2\in\{29,30,31,32\}, c_2^3  \in \{24,25,26,27\}, c_2^4 \in  \{19, 20,21,22\} $ and the base mobility matrix $\text{M}_4$ with $\text{M}_4[0,2] = \text{M}_4[2,0] = 0$.\\
 $B_t = 32 \times 10^3$ and $U_t^k = 10^4$.

Apart for the experiments in the main text, we provide more experiments in sections \ref{supp_sec:addition_exp_seir_K_3}, \ref{supp_sec:addition_exp_seir_K_4} and \ref{supp_sec:addition_exp_sepihr_K_4} to investigate the effect of changing the above parameters.

\section{Sampling Details} 
\subsection{ SEIR model }
In the sampling, we used the following upper and lower bounds on $\alpha,\beta \text{ and } \gamma$ respectively :\\
$0 \le \alpha \le 2$ \\ 
$ 0 \le \beta \le 1$ \\
$0 \le \gamma \le 1$

We ran the  sampling algorithm for a total of  $315$k iterations, from which we discard the first $15$k samples as burn-in \cite{fgpgm}.

\subsection{SEPIHR model}\label{supp_sec:sampling_sepihr}
We use $\alpha = 1.1, \beta = 0.08, \delta_1 = 0.01, \delta_2 = 0.002, \delta_3 =0.002, \gamma_1=0.1, \gamma_2 = 0.1 \text{ and } \gamma_3 =0.06 $ as the true parameter values. We simulate an SEPIHR model (eq:\eqref{eq:sepihr}) to get state values. We add  zero-mean Gaussian noise with $\sigma = 0.1$ to each of the simulated state values to generate our dataset. Using this data only for 15 days ($T=\{1,\dots,15\}$),  we first estimate the GP hyper-parameters for states using maximum-likelihood (see \eqref{eq:likelihood}).
We then run the Metropolis-Hastings MCMC sampling procedure to get our empirical posterior joint-distribution on $\alpha, \beta, \delta_1, \gamma_1 \text{ and } \gamma_2$. Note that as discussed earlier in section \ref{si_sec:sepihr_modelling}, we keep  $\delta_2 = 0.002, \delta_3 = 0.002 \text{ and } \gamma_3 = 0.06$ fixed, as these are known clinical parameters.
We use the following upper and lower bounds:\\
$0 \le \alpha \le 2$ \\
$ 0 \le \beta \le 1$\\
$0 \le \delta_1 \le 1$\\
$0 \le \gamma_1 \le 1$\\
$0 \le \gamma_2 \le 1$\\
After removing the burn-in samples, for the remaining  $3\times 10^5$ samples, we plot the marginal distributions along with their mean and mode in fig. \ref{fig:mcmc_dist_sepihr}. Mode is a more representative statistic here as the distributions are not symmetric due to non-negative support.

\begin{figure}[h]
\vspace{-0.1cm}
\label{fig:test}
\centering
\begin{subfigure}{\textwidth}
  \centering
  \includegraphics[width=\linewidth]{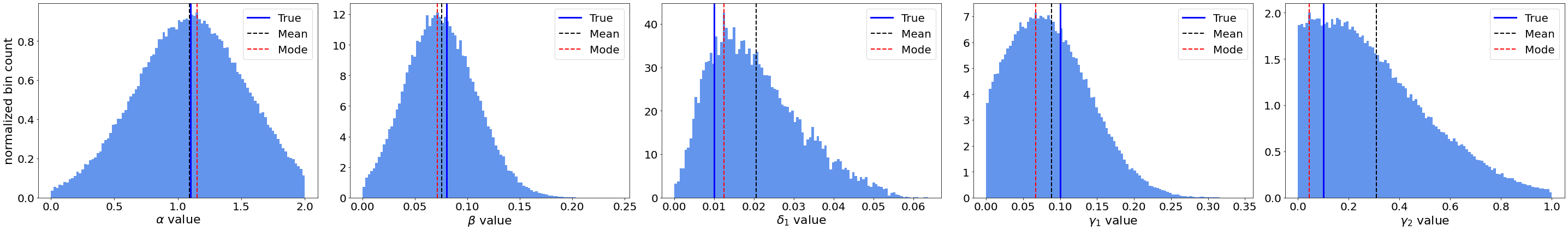}
\end{subfigure} 
\vspace{-0.6cm}
\caption{Empirical distribution after $3\times 10^5$ samples.}\label{fig:mcmc_dist_sepihr}
\end{figure}
\newpage

\section{Discussion on the resulting optimal vaccine policies}
\vspace{-0.3cm}
In this section we discuss the resulting optimal vaccine policies  i.e. nominal and stochastic policies for the three set of experiments in the main text.
\vspace{-0.2cm}
\subsection{SEIR model}
\begin{figure}[h!]
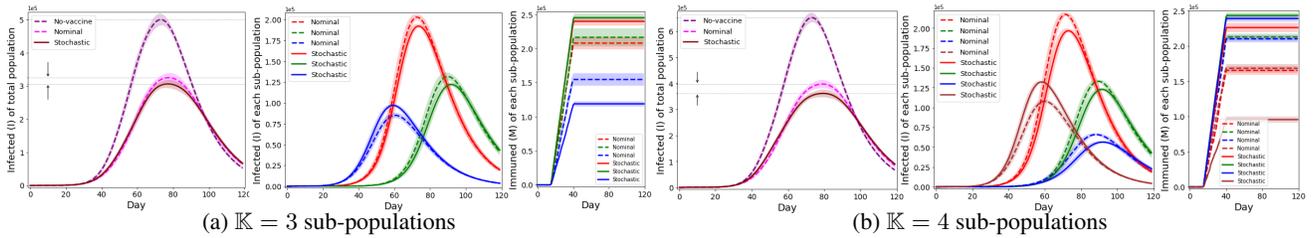

\begin{adjustwidth}{-0.4cm}{-0.4cm}
\centering
\captionsetup[subfigure]{aboveskip=-1pt,belowskip=-1pt}
\begin{subfigure}{.5\linewidth}
  \centering
  \includegraphics[width=\linewidth]{img/icml_figs_SEIR_pop_3_exps/rank_0_all.png}
  \caption{$\mathbb{K}=3$ sub-populations}
  \label{supp_fig:seirm_pop3}
\end{subfigure}%
\begin{subfigure}{.5\linewidth}
  \centering
  \includegraphics[width=\linewidth]{img/icml_figs_SEIR_pop_4_exps/rank_0_all.png}
  \caption{$\mathbb{K}=4$ sub-populations}
  \label{supp_fig:seirm_pop4}
\end{subfigure}
\vspace{-0.5cm}
\caption{\textbf{SEIR:} Evaluation of different vaccine policies i.e. no-vaccine $\mathcal{V}_{\phi}$, nominal $\mathcal{V}_{\mathcal{N}}$ and stochastic $\mathcal{V}_{\mathcal{S}}$.}
\label{supp_fig:seirm_exps}
\end{adjustwidth}
\end{figure}

For $\mathbb{K}=3$, we note from fig. \ref{supp_fig:seirm_pop3}, that under the nominal policy, least number of vaccines are allocated to the blue sub-population even though the onset of the pandemic is first in this sub-population. Interestingly,  even though the red sub-population is larger than the green sub-population, however, slightly more vaccines are allocated to the green sub-population in comparison to the red sub-population.
Under the stochastic policy, in the face of uncertainty, even fewer vaccines are allocated to the blue sub-population and proportionally more vaccine are allocated to red and blue sub-populations.\\ 
For $\mathbb{K}=4$, we note from fig. \ref{supp_fig:seirm_pop4}, that under the nominal policy, both the brown and the red sub-populations are allocated similar number of total vaccines despite the difference in their population size and the onset of the pandemic. Interestingly, even though the blue sub-population size is much smaller, however, it is allocated larger number of total vaccines in comparison to both the red and brown sub-populations.
Under the stochastic policy, the brown sub-population receives even fewer vaccines and proportionally more vaccines are allocated to the other three sub-populations.

\subsection{SEPIHR model}
\begin{figure}[h!]
\vspace{-0.3cm}
\centering
\begin{adjustwidth}{-0.3cm}{-0cm}
\begin{subfigure}{0.225\linewidth}
  \centering
  \includegraphics[width=0.98\linewidth]{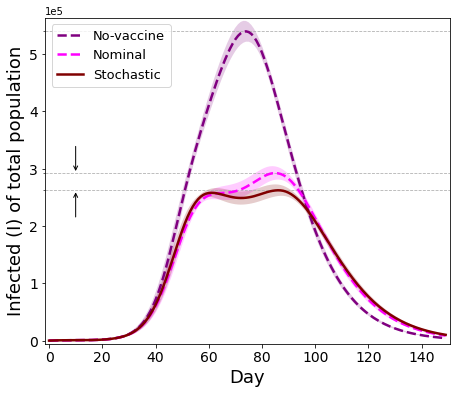}
\end{subfigure}%
\begin{subfigure}{0.225\linewidth}
  \centering
  \includegraphics[width=\linewidth]{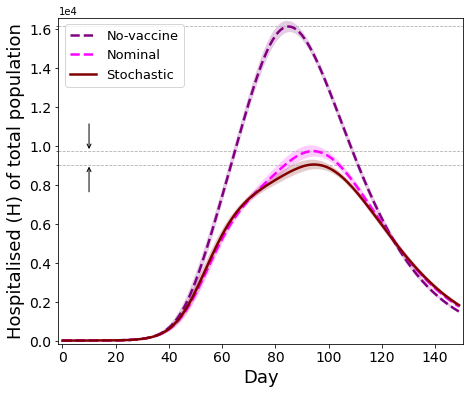}
\end{subfigure}%
\begin{subfigure}{0.225\linewidth}
  \centering
  \includegraphics[width=\linewidth]{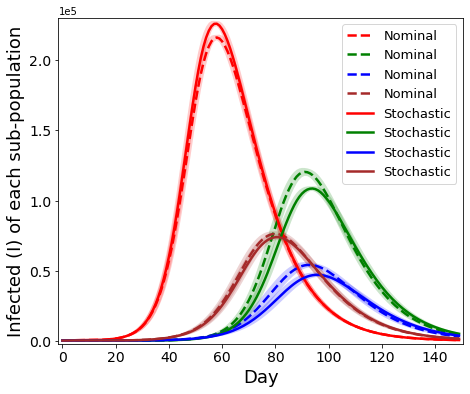}
\end{subfigure}%
\begin{subfigure}{0.225\linewidth}
  \centering
  \raisebox{0.0cm}{\includegraphics[width=0.96\linewidth]{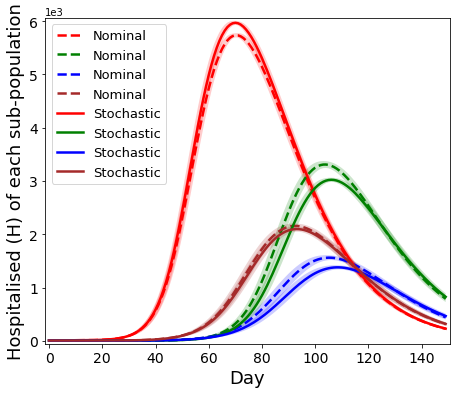}}
 \end{subfigure}%
 \begin{subfigure}{0.13\linewidth}
   \centering
  \includegraphics[width = \linewidth]{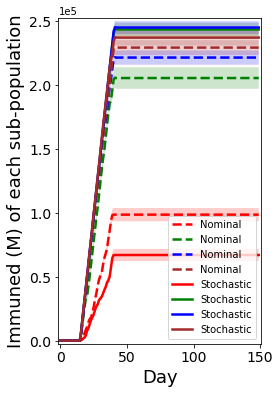}
 \end{subfigure}
\caption{\textbf{SEPIHR:} Evaluation of different policies: $\mathcal{V}_{\phi}$, $\mathcal{V}_{\mathcal{N}}$ and $\mathcal{V}_{\mathcal{S}}$ with $\mathbb{K}=4$ sub-populations.} \label{supp_fig:sepihrm_exps}
\end{adjustwidth}
\end{figure}

We note from fig. \ref{supp_fig:sepihrm_exps}, that under the nominal policy, least number of vaccines are allocated to the red sub-population regardless of its large size and the fact that the onset of the pandemic is first in this sub-population. This at first glance seems to be counter-intuitive, however can be explained by the fact that we aim to minimize the peak hospitalizations of the \emph{total} population which occurs around 82nd day.
Further, despite the  the brown and the blue sub-populations being smallest in size, receive larger number of vaccines.
Under the stochastic policy even fewer vaccines are allocated to red sub-population and correspondingly more vaccines are allocate to the other three sub-populations.
Importantly, also note that due to vaccination (both in the nominal and stochastic case) the peak of both the infections and hospitalisations is \textbf{delayed} in comparison to no-vaccination case,  thus further alleviating the burden of the medical personnel.

\section{Additional experiments for SEIR model with $\mathbb{K}=3$ sub-populations} \label{supp_sec:addition_exp_seir_K_3}
Evaluation of different vaccination policies i.e. no-vaccine policy $\mathcal{V}_{\phi}$, nominal policy $\mathcal{V}_{\mathcal{N}}$ and stochastic policy $\mathcal{V}_{\mathcal{S}}$ under different scenarios of $c_2^1,c_2^2, c_2^3 $ and mobility matrix M.

\begin{figure*}[h!]
\vspace{-0.5cm}
\centering
\begin{subfigure}{.33\textwidth}
  \centering
  \includegraphics[width=1.06\linewidth]{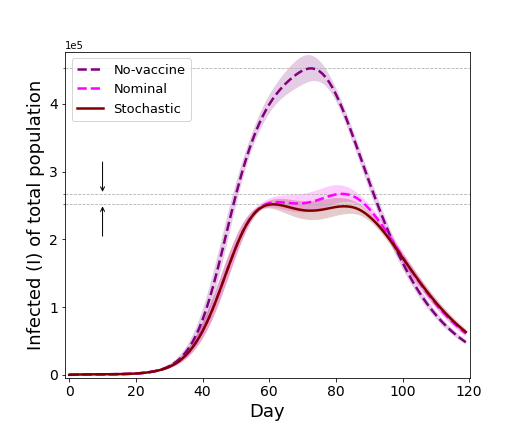}
  \caption{}
  \label{fig:tot_I}
\end{subfigure}%
\begin{subfigure}{0.33\textwidth}
  \centering
  \includegraphics[width=\linewidth]{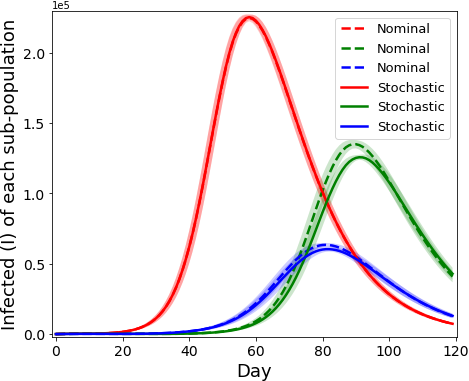}
  \caption{}
  \label{fig:sub_I}
 \end{subfigure}
 \begin{subfigure}{0.33\textwidth}
   \centering
  \includegraphics[width=0.985\linewidth]{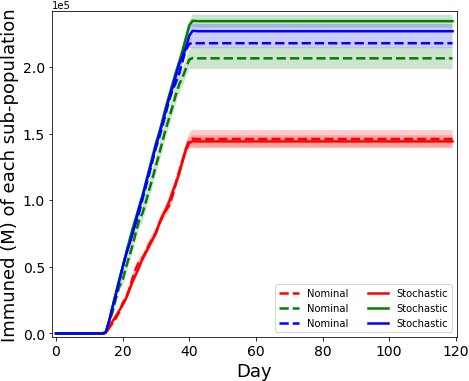}
   \caption{}
   \label{fig:sub_M}
 \end{subfigure}
\vspace{-0.1cm}
\caption{Peak infections (I) for the three policies (i.e. $\mathcal{V}_{\phi}$, $\mathcal{V}_{\mathcal{N}}$ and $\mathcal{V}_{\mathcal{S}}$) are around 453189, 264453 and 255233 respectively. $\mathcal{V}_{\mathcal{S}}$ provides a reduction of \textbf{3.486\%} in peak infections (I) over $\mathcal{V}_{\mathcal{N}}$.\\
Setup: $ c_2^1 \in \{ 9.0,10,11.0,12.0 \} , c_2^2 \in \{ 29.0,30,31.0,32.0 \} , c_2^3 \in \{ 19.0,20,21.0,22.0 \} $ and mobility matrix $\text{M}_3$  with $\text{M}_3[0,1] = 0$ and $\text{M}_3[1,0] = 0$.
 }
\label{fig:test}
\end{figure*}

\vspace{-0.55cm}

\begin{figure*}[h!]
\centering
\begin{subfigure}{.33\textwidth}
  \centering
  \includegraphics[width=1.06\linewidth]{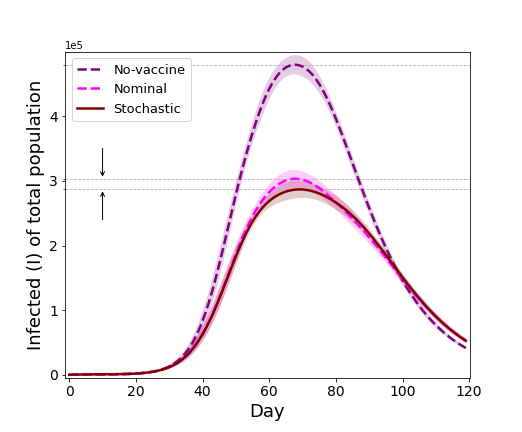}
  \caption{}
  \label{fig:tot_I}
\end{subfigure}%
\begin{subfigure}{0.33\textwidth}
  \centering
  \includegraphics[width=\linewidth]{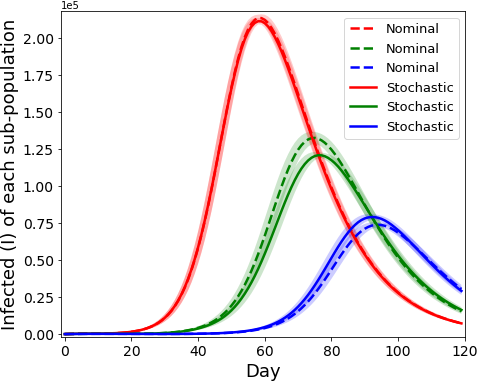}
  \caption{}
  \label{fig:sub_I}
 \end{subfigure}
 \begin{subfigure}{0.33\textwidth}
   \centering
  \includegraphics[width=0.985\linewidth]{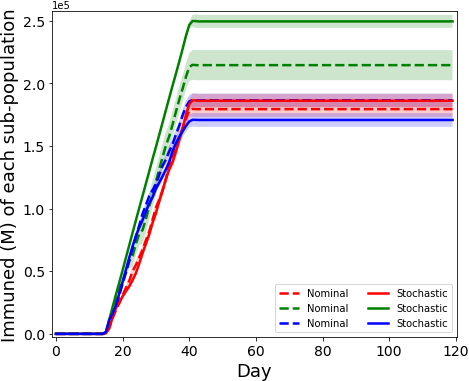}
   \caption{}
   \label{fig:sub_M}
 \end{subfigure}
\vspace{-0.1cm}
\caption{Peak infections (I) for the three policies (i.e. $\mathcal{V}_{\phi}$, $\mathcal{V}_{\mathcal{N}}$ and $\mathcal{V}_{\mathcal{S}}$) are around 481463, 298401 and 288031 respectively. $\mathcal{V}_{\mathcal{S}}$ provides a reduction of \textbf{3.475\%} in peak infections (I) over $\mathcal{V}_{\mathcal{N}}$.\\
Setup: $ c_2^1 \in \{ 9.0,10,11.0,12.0 \} , c_2^2 \in \{ 19.0,20,21.0,22.0 \} , c_2^3 \in \{ 29.0,30,31.0,32.0 \} $ and mobility matrix $\text{M}_3$  with $\text{M}_3[0,1] = 0$ and $\text{M}_3[1,0] = 0$.
 }
\label{fig:test}
\vspace{-0.65cm}
\end{figure*}


\begin{figure*}[h!]
\centering
\begin{subfigure}{.33\textwidth}
  \centering
  \includegraphics[width=1.06\linewidth]{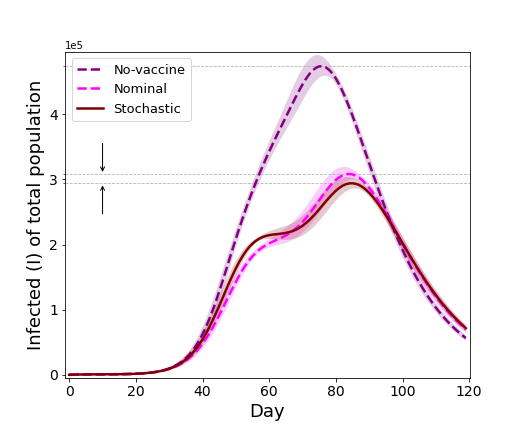}
  \caption{}
  \label{fig:tot_I}
\end{subfigure}%
\begin{subfigure}{0.33\textwidth}
  \centering
  \includegraphics[width=\linewidth]{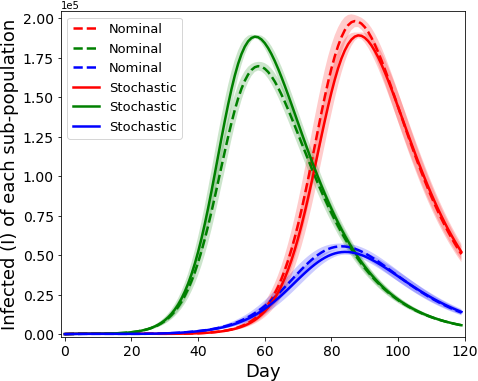}
  \caption{}
  \label{fig:sub_I}
 \end{subfigure}
 \begin{subfigure}{0.33\textwidth}
   \centering
  \includegraphics[width=0.985\linewidth]{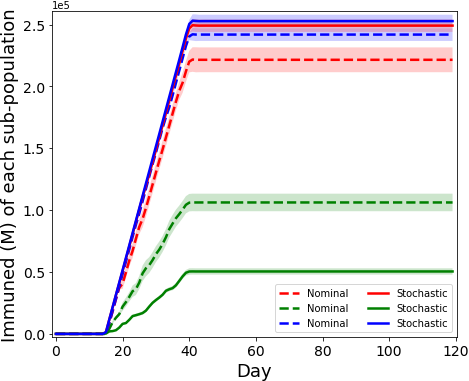}
   \caption{}
   \label{fig:sub_M}
 \end{subfigure}
\vspace{-0.1cm}
\caption{Peak infections (I) for the three policies (i.e. $\mathcal{V}_{\phi}$, $\mathcal{V}_{\mathcal{N}}$ and $\mathcal{V}_{\mathcal{S}}$) are around 475956, 305852 and 295492 respectively. $\mathcal{V}_{\mathcal{S}}$ provides a reduction of \textbf{3.387\%} in peak infections (I) over $\mathcal{V}_{\mathcal{N}}$.\\
Setup: $ c_2^1 \in \{ 29.0,30,31.0,32.0 \} , c_2^2 \in \{ 9.0,10,11.0,12.0 \} , c_2^3 \in \{ 19.0,20,21.0,22.0 \} $ and mobility matrix $\text{M}_3$  with $\text{M}_3[0,2] = 0$ and $\text{M}_3[2,0] = 0$.
}
\label{fig:test}
\end{figure*}

\clearpage
\section{Additional experiments for SEIR model with $\mathbb{K}=4$ sub-populations}\label{supp_sec:addition_exp_seir_K_4}
Evaluation of different vaccination policies i.e. no-vaccine policy $\mathcal{V}_{\phi}$, nominal policy $\mathcal{V}_{\mathcal{N}}$ and stochastic policy $\mathcal{V}_{\mathcal{S}}$ under different scenarios of $c_2^1,c_2^2, c_2^3,c_2^4 $ and mobility matrix M. 

\begin{figure*}[h!]
\vspace{-0.25cm}
\centering
\begin{subfigure}{.33\textwidth}
  \centering
  \includegraphics[width=0.95\linewidth]{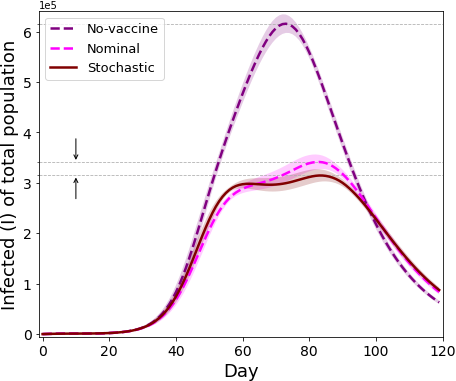}
  \caption{}
  \label{fig:rank_1_tot_I}
\end{subfigure}%
\begin{subfigure}{0.33\textwidth}
  \centering
  \includegraphics[width=\linewidth]{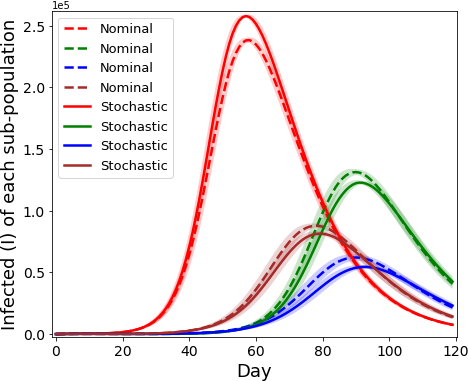}
  \caption{}
  \label{fig:rank_1_sub_I}
 \end{subfigure}
 \begin{subfigure}{0.33\textwidth}
   \centering
  \includegraphics[width=\linewidth]{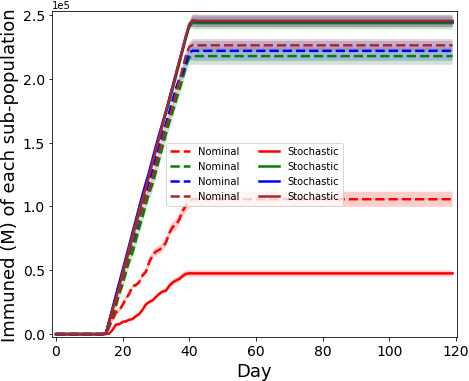}
   \caption{}
   \label{fig:rank_1_sub_M}
 \end{subfigure}
\vspace{-0.1cm}
\caption{Peak infections (I) for the three policies (i.e. $\mathcal{V}_{\phi}$, $\mathcal{V}_{\mathcal{N}}$ and $\mathcal{V}_{\mathcal{S}}$) are around 616736, 342776 and 316411 respectively. $\mathcal{V}_{\mathcal{S}}$ provides a reduction of \textbf{7.692\%} in peak infections (I) over $\mathcal{V}_{\mathcal{N}}$. \\
Setup: $ c_2^1 \in \{ 9.0,10,11.0,12.0 \} , c_2^2 \in \{ 29.0,30,31.0,32.0 \} , c_2^3 \in \{ 24.0,25,26.0,27.0 \} , c_2^4 \in \{ 19.0,20,21.0,22.0 \} $ and mobility matrix $\text{M}_4$  with $\text{M}_4[1,3] = 0$ and $\text{M}_4[3,1] = 0$. }
\label{fig:rank_1_seir_k_4}
\end{figure*}

\vspace{-0.25cm}

\begin{figure*}[h!]
\centering
\begin{subfigure}{.33\textwidth}
  \centering
  \includegraphics[width=0.95\linewidth]{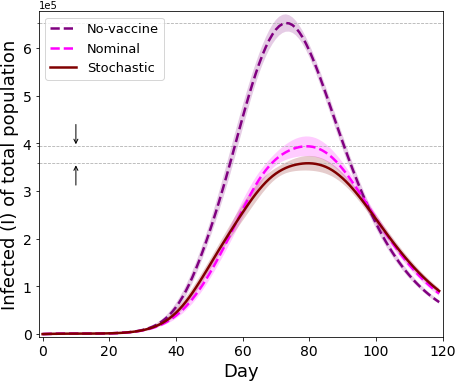}
  \caption{}
  \label{fig:rank_5_tot_I}
\end{subfigure}%
\begin{subfigure}{0.33\textwidth}
  \centering
  \includegraphics[width=\linewidth]{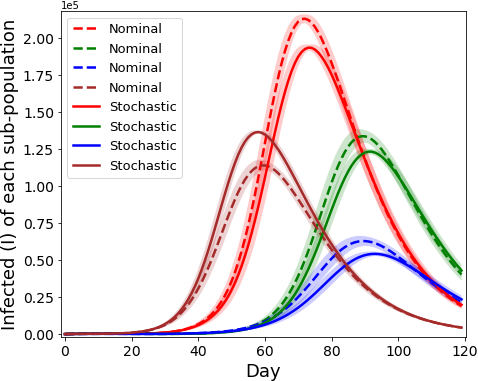}
  \caption{}
  \label{fig:rank_5_sub_I}
 \end{subfigure}
 \begin{subfigure}{0.33\textwidth}
   \centering
  \includegraphics[width=\linewidth]{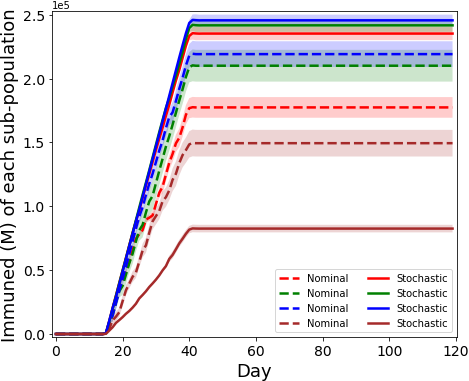}
   \caption{}
   \label{fig:rank_5_sub_M}
 \end{subfigure}
\vspace{-0.1cm}
\caption{Peak infections (I) for the three policies (i.e. $\mathcal{V}_{\phi}$, $\mathcal{V}_{\mathcal{N}}$ and $\mathcal{V}_{\mathcal{S}}$) are around 653693, 383850 and 360789 respectively. $\mathcal{V}_{\mathcal{S}}$ provides a reduction of \textbf{6.008\%} in peak infections (I) over $\mathcal{V}_{\mathcal{N}}$.\\
Setup: $ c_2^1 \in \{ 19.0,20,21.0,22.0 \} , c_2^2 \in \{ 29.0,30,31.0,32.0 \} , c_2^3 \in \{ 24.0,25,26.0,27.0 \} , c_2^4 \in \{ 9.0,10,11.0,12.0 \} $ and mobility matrix $\text{M}_4$  with $\text{M}_4[0,2] = 0$ and $\text{M}_4[2,0] = 0$.
  }
\label{fig:rank_5_seir_k_4}
\vspace{-0.25cm}
\end{figure*}
\begin{figure*}[ht!]
\centering
\begin{subfigure}{.33\textwidth}
  \centering
  \includegraphics[width=0.95\linewidth]{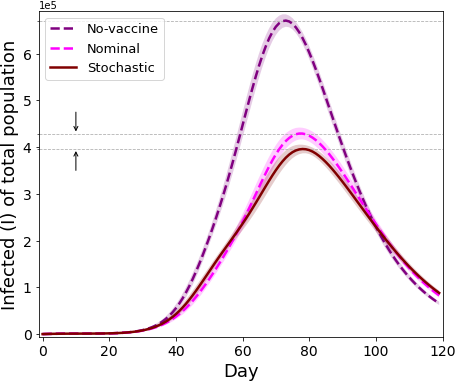}
  \caption{}
  \label{fig:rank_6_tot_I}
\end{subfigure}%
\begin{subfigure}{0.33\textwidth}
  \centering
  \includegraphics[width=\linewidth]{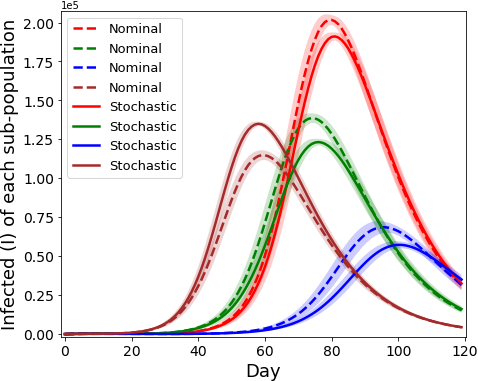}
  \caption{}
  \label{fig:rank_6_sub_I}
 \end{subfigure}
 \begin{subfigure}{0.33\textwidth}
   \centering
  \includegraphics[width=\linewidth]{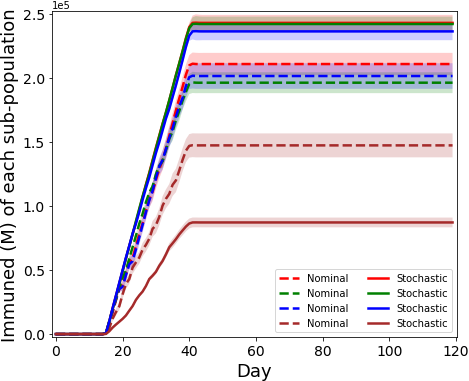}
   \caption{}
   \label{fig:rank_6_sub_M}
 \end{subfigure}
\vspace{-0.1cm}
\caption{Peak infections (I) for the three policies (i.e. $\mathcal{V}_{\phi}$, $\mathcal{V}_{\mathcal{N}}$ and $\mathcal{V}_{\mathcal{S}}$) are around 671850, 421908 and 397524 respectively. $\mathcal{V}_{\mathcal{S}}$ provides a reduction of \textbf{5.779\%} in peak infections (I) over $\mathcal{V}_{\mathcal{N}}$.\\
Setup: $ c_2^1 \in \{ 24.0,25,26.0,27.0 \} , c_2^2 \in \{ 19.0,20,21.0,22.0 \} , c_2^3 \in \{ 29.0,30,31.0,32.0 \} , c_2^4 \in \{ 9.0,10,11.0,12.0 \} $ and mobility matrix $\text{M}_4$ .
  }
\label{fig:rank_6_seir_k_4}
\end{figure*}

\clearpage
\section{SEPIHR model details \& Nominal and Stochastic optimization problem formulations}\label{si_sec:sepihr_modelling}
\vspace{-0.1cm}
We mathematically describe the ODEs governing the evolution of each of the states in the SEPIHR disease transmission model shown in fig. \ref{fig:si_sepihr} in equations \eqref{eq:sepihr}.
\begin{adjustwidth}{0cm}{-1cm}
\begin{minipage}{0.5\linewidth}
\begin{subequations}\label{eq:sepihr}
\begin{alignat}{10}
&\frac{dS(t)}{dt}   &&\coloneqq \dot{S}(t) =  - \frac{\alpha}{N} S(t) I(t)          \label{eq:sepihr_s}\\
&\frac{dE(t)}{dt}   &&\coloneqq \dot{E}(t) = \frac{\alpha}{N} S(t) I(t) -  (\beta + \delta_1)E(t)  \label{eq:sepihr_e}\\
& \frac{dP(t)}{dt} &&\coloneqq \dot{P}(t) =  \delta_1 E(t) -(\delta_2 + \gamma_2)P(t)  \label{eq:sepihr_p}\\
&\frac{dI(t)}{dt}   &&\coloneqq \dot{I}(t) = \beta E(t) - (\gamma_1 +\delta_3)  I(t)              \label{eq:sepihr_i}\\
& \frac{dH(t)}{dt} &&\coloneqq \dot{H}(t) = \delta_2P(t) + \delta_3I(t)  - \gamma_3H(t) \label{eq:sepihr_h}\\
&\frac{dR(t)}{dt}   &&\coloneqq \dot{R}(t) = \gamma_1 I(t) +  \gamma_2 P(t) +\gamma_3 H(t)   \phantom{MMM}                     \label{eq:sepihr_r}
\end{alignat}
\end{subequations}
\end{minipage}
\begin{minipage}{0.5\linewidth}
\centering 
\includegraphics[scale = 0.6]{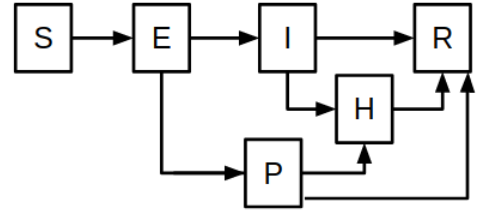}
 \captionof{figure}{SEPIHR model.}\label{fig:si_sepihr}
\end{minipage}
\end{adjustwidth}
\vspace{0.25cm}
The SEPIHR model \eqref{eq:sepihr} has two set of parameters, where in the first set we have $\alpha, \beta, \delta_1, \gamma_1, \gamma_2$ and in the second set we have $\delta_2, \delta_3,  \gamma_3$. Parameters in the second set, also referred to as clinical parameters, i.e. $\delta_2, \delta_3$ and $\gamma_3$ are the rates at which a quarantined person is hospitalised, an infected person is hospitalised and the rate at which a hospitalised person leaves the hospital respectively. These rates (i.e. second set of parameters) have standard values \cite{bertsimas_vaccine}, for instance a person on an average would spend 15-16 days in the hospital, therefore $\gamma_3 = 1/16 \approx 0.06$. Importantly, accurately estimating these clinical parameters require data which is typically only available to hospitals and medical personnel. Therefore we prefer to use the standard known values for these parameters instead of estimating them from daily counts data. We use $\delta_2 = 0.002, \delta_3 =0.002, \text{ and } \gamma_3 =0.06 $ as the true parameter values. 
Additionally, this also helps to keep the sampling space small (i.e. 5 dimensional), therefore MCMC sampling remains computationally tractable as we only estimate the posterior joint distribution over the first set of parameters, i.e. $\alpha, \beta, \delta_1, \gamma_1 \text{ and } \gamma_2$. For sampling details, refer to section \ref{supp_sec:sampling_sepihr}.

The spread of disease in each of the $k$-th sub-population, where $k \in \{1,\dots,\mathbb{K}\}$, is modeled using a separate SEPIHR model. 
To account for the vaccinated individuals, the SEPIHR model in fig. \ref{fig:si_sepihr} is updated with a new compartment (denoted by M) to represent the immune population and the updated model shown in fig. \ref{fig:si_sepihrm}, is denoted as the SEPIHRM model. Let $V_k(t)$ represent the number of people vaccinated at time $t$ in the $k$-th sub-population  and $\eta$ be efficacy of the vaccine, then the ODEs corresponding to the SEIRM model of the $k$-th sub-population  are given as:
\begin{adjustwidth}{0cm}{-1cm}
\begin{minipage}{0.5\linewidth}
\begin{subequations}\label{eq:seirm}
\begin{alignat}{10}
&\frac{dS_k(t)}{dt}   &&\coloneqq \dot{S}_k(t) =  - \frac{\alpha}{N} ( S_k(t) - \eta V_k(t)) I_k(t)          \label{eq:sepihrm_s}\\
&\frac{dE_k(t)}{dt}   &&\coloneqq \dot{E}_k(t) = \frac{\alpha}{N} ( S_k(t)- \eta V_k(t)) I_k(t) -  (\beta + \delta_1)E_k(t)  \label{eq:sepihrm_e}\\
& \frac{dP_k(t)}{dt} &&\coloneqq \dot{P}_k(t) =  \delta_1 E_k(t) -(\delta_2 + \gamma_2)P_k(t)  \label{eq:sepihrm_p}\\
&\frac{dI_k(t)}{dt}   &&\coloneqq \dot{I}_k(t) = \beta E_k(t) - (\gamma_1 +\delta_3)  I_k(t)              \label{eq:sepihrm_i}\\
& \frac{dH_k(t)}{dt} &&\coloneqq \dot{H}_k(t) =\delta_2P_k(t)+ \delta_3I_k(t)   - \gamma_3H_k(t)   \label{eq:sepihrm_h}\\
&\frac{dR_k(t)}{dt}   &&\coloneqq \dot{R}_k(t) = \gamma_1 I_k(t) +\gamma_2 P_k(t)+ \gamma_3 H_k(t)    \label{eq:sepihrm_r}  \\
& \frac{dM_k(t)}{dt} && \coloneqq \dot{M}_k(t) = \eta V_k(t)\label{eq:sepihrm_m}
\end{alignat}
\end{subequations}
\end{minipage}
\begin{minipage}{0.5\linewidth}
\centering 
\includegraphics[scale = 0.6]{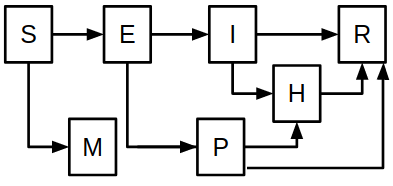}
 \captionof{figure}{SEPIHRM model.}\label{fig:si_sepihrm}
\end{minipage}
\end{adjustwidth}
\vspace{0.25cm}
We consider two important features of disease transmission. First, due to mobility there is contact between infected individuals of one sub-population with the susceptible individuals of another sub-population. Second, due to different levels of mobility between different sub-populations,
onset of the pandemic in each of the sub-populations generally vary. Both of these are accounted by updating the states $S_k$ and $E_k$ as follows:
\begin{adjustwidth}{-0.4cm}{-0cm}
\begin{minipage}{0.73\textwidth}
\begin{subequations}\label{eq:seirm}
\begin{small}
\begin{alignat}{10}
&\frac{dS_k(t)}{dt}   &&\coloneqq \dot{S}_k(t) =  - \frac{u_k(t)  \alpha}{N} \big( S_k(t) - \eta V_k(t)\big)\big(\sum_{r=1}^{K} {\lambda}_{r}^{k}I_r(t)  \big)          \label{eq:updated_sepihrm_s}\\
&\frac{dE_k(t)}{dt}   &&\coloneqq \dot{E}_k(t) = \;\;\; \frac{u_k(t) \alpha}{N} \big( S_k(t)- \eta V_k(t)\big)\big(\sum_{r=1}^{K} {\lambda}_{r}^{k}I_r(t)  \big) -  (\beta + \delta_1)E_k(t)  \label{eq:updated_sepihrm_e}
\end{alignat}
\end{small}
\end{subequations}
\end{minipage}
\begin{minipage}{0.25\textwidth}
\end{minipage}
\end{adjustwidth}

where 
${\lambda}_{r}^{k}$ denotes the mobility levels from sub-population $r$ to sub-population $k$ and $u_k(t)$ corresponds to a sigmoid function, $u_k(t) \coloneqq 1/(1 + e^{-c_1^k(t-c_2^k)} )$ with parameters $c_1^k \text{ and } c_2^k$.
In particular, $c_2^k$ controls the onset of the pandemic in the $k$-th sub-population, therefore  we also account for uncertainty in $c_2^k \; \forall\; k \in \mathcal{K}$, by appropriately extending the scenario set $\Omega$.
We can now write the corresponding \emph{nominal (or non-stochastic)} optimization problem (denoted $\mathcal{NF}$) for vaccine allocation as follows:
\begin{adjustwidth}{-0.1cm}{-0cm}
\begin{minipage}{0.9\linewidth}
\begin{subequations}
\begin{alignat}{10}
\mathcal{NF}: & \min_{V} \;\; \mathcal{H} \phantom{MMMMMMMMMMMMMMMMM} \text{(Nominal Formulation)} \label{eq:nominal2_obj} \\
\text{s.t.} \;\;  & \notag \\
&  \{\text{\eqref{eq:updated_sepihrm_s},\eqref{eq:updated_sepihrm_e},\eqref{eq:sepihrm_p},\eqref{eq:sepihrm_i},\eqref{eq:sepihrm_h},\eqref{eq:sepihrm_r},\eqref{eq:sepihrm_m}}   \}  \phantom{MMMMM}  \forall \;  k \in \mathcal{K}, t \in \mathcal{T} \label{eq:nominal2_odes} \\
&   \sum_{k=1}^{\mathbb{K}} H_k(t)  = \tilde{H}(t) \phantom{MMMMMMMMMMMMMM} \forall \; \; t \in \mathcal{T}  \label{eq:nominal2_Isum}\\
&   \tilde{H}(t)  \le \mathcal{H} \phantom{MMMMMMMMMMMMMMMMM} \forall  \;\; t \in \mathcal{T} \label{eq:nominal2_Imax}\\
&  \sum_{k=1}^{\mathbb{K}} V_k(t)  \le B_{t} \phantom{MMMMMMMMMMMMMMm} \forall \;\;  t \in \{t_{s},\dots,t_{l}\} \label{eq:nominal2_total_B}\\
& 0 \le V_{k}(t) \le U_{t}^k  \phantom{MMMMMMMMMMMMMmm} \forall \;\; k \in \mathcal{K}, t \in \{t_{s},\dots,t_{l} \} \label{eq:nominal2_ind_B} \\
&  V_{k}(t) = 0  \phantom{MMMMMMMMMMMMMMMM} \;\;\; \forall \;\; k \in \mathcal{K}, t \in \mathcal{T} \setminus  \{t_s,\dots, t_l \} \label{eq:nominal2_v_0}
\end{alignat}
\end{subequations}
\end{minipage}
\end{adjustwidth}

where $B_t$ denotes the daily total vaccination budget and $U_t^k$ denotes the daily vaccination budget of the $k$-th sub-population. Note that in the above formulations, we are minimizing peak hospitalizations. We next provide the \textbf{uncertainty-informed}, i.e. stochastic counterpart (denoted $\mathcal{SF}$) of the above nominal problem $\mathcal{NF}$. 

\begin{subequations}
\begin{alignat}{10}
\hypertarget{SF}{\mathcal{SF}}: & \min_{V} \quad \sum_{\omega \in \Omega } p_{\omega}\mathcal{H}_{\omega} \phantom{MMMMMMMMMMMMMm} \text{(Stochastic Formulation)}  \label{eq:stoc2_obj} \\
\text{s.t.} \;\; & \notag \\
& \begin{rcases} \dot{S}_k^{\omega}(t) =  - \frac{u_k^{\omega}(t)  \alpha^{\omega}}{N} \big( S_k^{\omega}(t) - \eta V_k(t)\big)\big(\sum_{r=1}^{\mathbb{K}} {\lambda}_{r}^{k}I_r^{\omega}(t)  \big)        \\
 \dot{E}_k^{\omega}(t) = \;\;\; \frac{u_k^{\omega}(t) \alpha^{\omega}}{N} \big( S_k^{\omega}(t)- \eta V_k^{\omega}(t)\big)\big(\sum_{r=1}^{\mathbb{K}} {\lambda}_{r}^{k}I_r^{\omega}(t)  \big) -  (\beta^{\omega} + \delta_1^{\omega})E_k^{\omega}(t) \\
\dot{P}_k^{\omega}(t) =  \delta_1^{\omega} E_k^{\omega}(t) -(\delta_3 + \gamma_3)P_k^{\omega}(t)  \\
 \dot{I}_k^{\omega}(t) = \beta^{\omega} E_k^{\omega}(t) - (\gamma_1^{\omega} +\delta_2)  I_k^{\omega}(t)             \\
 \dot{H}_k^{\omega}(t) = \delta_2I_k^{\omega}(t) + \delta_3 Q_k^{\omega}(t) - \gamma_2^{\omega}H_k^{\omega}(t) \\
 \dot{R}_k^{\omega}(t) = \gamma_1^{\omega} I_k^{\omega}(t) + \gamma_2^{\omega} H_k^{\omega}(t) + \gamma_3 Q_k^{\omega}(t)                         \\
\dot{M}_k^{\omega}(t) = \eta V_k(t)\\
 \sum_{k=1}^{K} H_k^{\omega}(t)  = \tilde{H}^{\omega}(t) \phantom{M} \\ 
 \tilde{H}^{\omega}(t)  \le \mathcal{H}_{\omega}  \end{rcases} \forall k \in \mathcal{K}, t \in \mathcal{T}, \omega \in \Omega  \\
 & \{\text{\eqref{eq:nominal2_total_B}, \eqref{eq:nominal2_ind_B}, \eqref{eq:nominal2_v_0}} \}   \label{eq:stoc2_budget_v}
\end{alignat}
\end{subequations}

\newpage
\section{Additional experiments for SEPIHR model with $\mathbb{K}=4$ sub populations}\label{supp_sec:addition_exp_sepihr_K_4}
\vspace{-0.1cm}
\begin{figure}[h!]
\centering
\begin{subfigure}{.4\textwidth}
  \centering
  \includegraphics[width=0.99\linewidth]{img/icml_figs_SEQIHR_pop_4_exps/rank_1_infected_total_I.png}
  \label{fig:rank_1_tot_I}
\end{subfigure}%
\hspace{1cm}
\begin{subfigure}{.4\textwidth}
  \centering
  \includegraphics[width=\linewidth]{img/icml_figs_SEQIHR_pop_4_exps/rank_1_infected_total_H.png}
  \label{fig:rank_1_tot_H}
\end{subfigure}%
\vspace{-0.5cm}
\begin{subfigure}{.4\textwidth}
  \centering
  \includegraphics[width=\linewidth]{img/icml_figs_SEQIHR_pop_4_exps/rank_1_pop_individual_I.png}
  \label{fig:rank_1_ind_I}
\end{subfigure}%
\hspace{1cm}
\raisebox{0.45cm}{
\begin{subfigure}{0.4\textwidth}
  \centering
  \includegraphics[width=0.96\linewidth]{img/icml_figs_SEQIHR_pop_4_exps/rank_1_pop_individual_H.png}
  \label{fig:rank_1_ind_H}
 \end{subfigure}}%
 \vspace{-0.5cm}
 \begin{subfigure}{0.4\textwidth}
   \centering
  \includegraphics[width=0.99\linewidth]{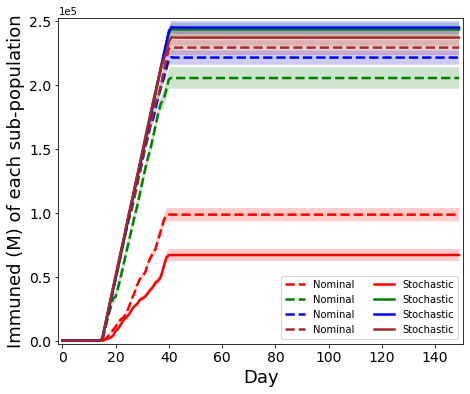}
   \label{fig:rank_1_ind_M}
 \end{subfigure}
 \vspace{-0.6cm}
\caption{Peak infections (I) for the three policies (i.e. $\mathcal{V}_{\phi}$, $\mathcal{V}_{\mathcal{N}}$ and $\mathcal{V}_{\mathcal{S}}$) are around 539451, 280526 and 262650 and the peak hospitalisations are around 16122, 9477 and 9061 respectively. $\mathcal{V}_{\mathcal{S}}$ provides a reduction of \textbf{6.37\%} in peak infections (I) over $\mathcal{V}_{\mathcal{N}}$ and a \textbf{4.4\%} reduction in peak hospitalisations (H).\\
Setup: $ c_2^1 \in \{ 9.0,10,11.0,12.0 \} , c_2^2 \in \{ 29.0,30,31.0,32.0 \} , c_2^3 \in \{ 24.0,25,26.0,27.0 \} , c_2^4 \in \{ 19.0,20,21.0,22.0 \} $ and mobility matrix $\text{M}_4$  with $\text{M}_4[0,2] = 0$ and $\text{M}_4[2,0] = 0$.
 }
\label{fig:rank_1_seqihr}
\end{figure}
\newpage
\begin{figure}[h!]
\centering
\begin{subfigure}{.425\textwidth}
  \centering
  \includegraphics[width=0.99\linewidth]{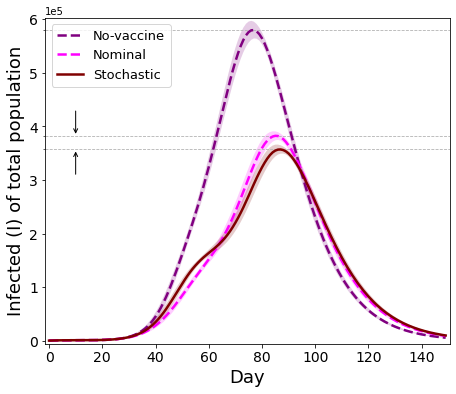}
  \label{fig:rank_12_tot_I}
\end{subfigure}%
\hspace{1cm}
\begin{subfigure}{.425\textwidth}
  \centering
  \includegraphics[width=\linewidth]{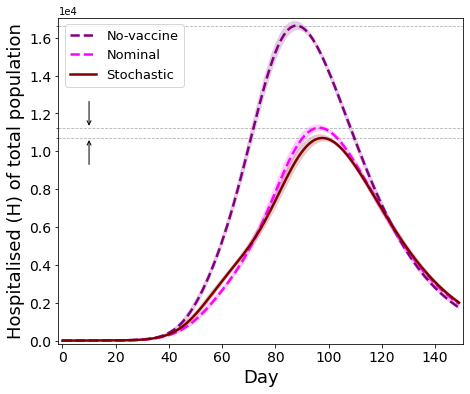}
  \label{fig:rank_12_tot_H}
\end{subfigure}%
\vspace{-0.5cm}
\begin{subfigure}{.425\textwidth}
  \centering
  \includegraphics[width=\linewidth]{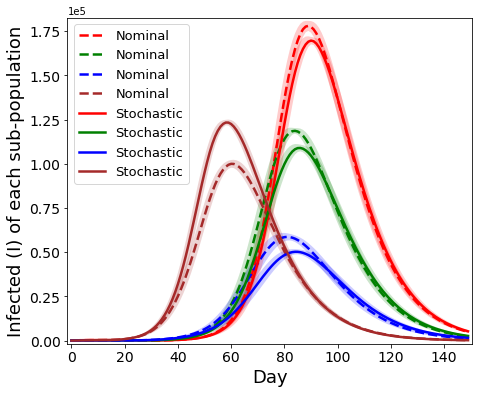}
  \label{fig:rank_12_ind_I}
\end{subfigure}%
 \hspace{0.9cm}
\raisebox{0.38cm}{
\begin{subfigure}{0.425\textwidth}
  \centering
  \includegraphics[width=0.96\linewidth]{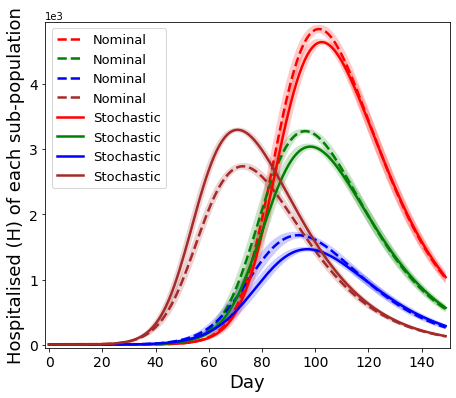}
  \label{fig:rank_12_ind_H}
 \end{subfigure}}%
 \vspace{-0.5cm}
 \begin{subfigure}{0.45\textwidth}
   \centering
  \includegraphics[width=0.99\linewidth]{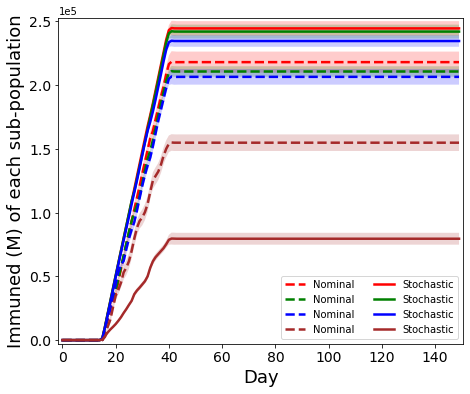}
   \label{fig:rank_12_ind_M}
 \end{subfigure}
\vspace{-0.5cm} 
\caption{Peak infections (I) for the three policies (i.e. $\mathcal{V}_{\phi}$, $\mathcal{V}_{\mathcal{N}}$ and $\mathcal{V}_{\mathcal{S}}$) are around 580009, 377703 and 357624 and the peak hospitalisations are around 16651, 11160 and 10725 respectively. $\mathcal{V}_{\mathcal{S}}$ provides a reduction of \textbf{5.32\%} in peak infections (I) over $\mathcal{V}_{\mathcal{N}}$ and a \textbf{3.9\%} reduction in peak hospitalisations (H).\\
Setup: $ c_2^1 \in \{ 29.0,30,31.0,32.0 \} , c_2^2 \in \{ 24.0,25,26.0,27.0 \} , c_2^3 \in \{ 19.0,20,21.0,22.0 \} , c_2^4 \in \{ 9.0,10,11.0,12.0 \} $ and mobility matrix $\text{M}_4$  with $\text{M}_4[0,2] = 0$ and $\text{M}_4[2,0] = 0$. }
\label{fig:rank_12_seqihr}
\end{figure}
\newpage
\begin{figure}[h!]
\centering
\begin{subfigure}{.425\textwidth}
  \centering
  \includegraphics[width=0.99\linewidth]{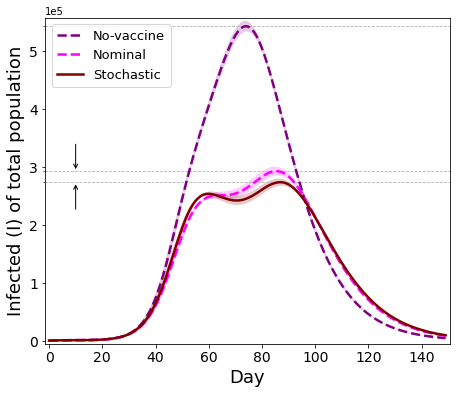}
  \label{fig:rank_13_tot_I}
\end{subfigure}%
\hspace{1cm}
\begin{subfigure}{.425\textwidth}
  \centering
  \includegraphics[width=\linewidth]{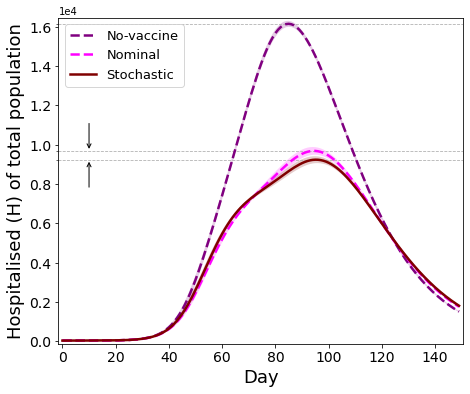}
  \label{fig:rank_13_tot_H}
\end{subfigure}%
\vspace{-0.5cm}
\begin{subfigure}{.425\textwidth}
  \centering
  \includegraphics[width=\linewidth]{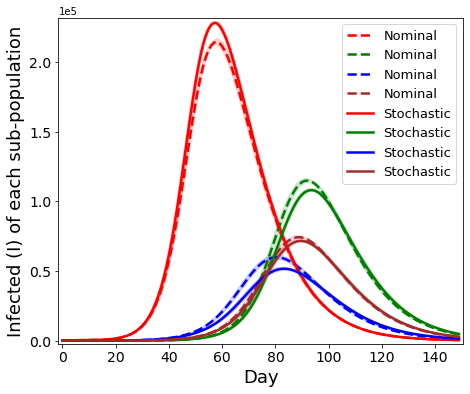}
  \label{fig:rank_13_ind_I}
\end{subfigure}%
\hspace{1cm}
\raisebox{0.45cm}{
\begin{subfigure}{0.425\textwidth}
  \centering
  \includegraphics[width=0.96\linewidth]{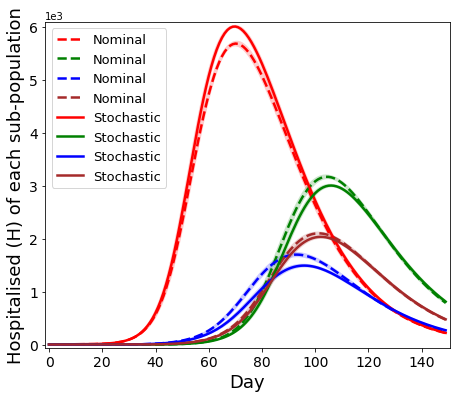}
  \label{fig:rank_13_ind_H}
 \end{subfigure}}
 \begin{subfigure}{0.425\textwidth}
   \centering
  \includegraphics[width=0.99\linewidth]{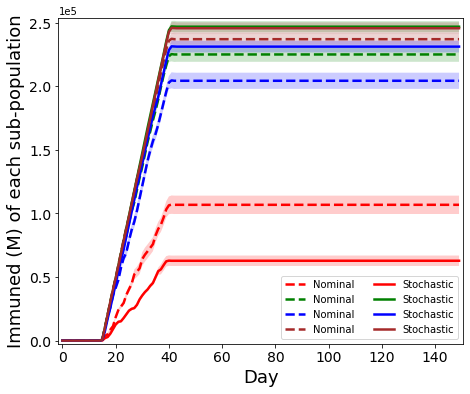}
   \label{fig:rank_13_ind_M}
 \end{subfigure}
 \vspace{-0.5cm}
\caption{Peak infections (I) for the three policies (i.e. $\mathcal{V}_{\phi}$, $\mathcal{V}_{\mathcal{N}}$ and $\mathcal{V}_{\mathcal{S}}$) are around 543114, 289457 and 274104 and the peak hospitalisations are around 16152, 9625 and 9244 respectively. $\mathcal{V}_{\mathcal{S}}$ provides a reduction of \textbf{5.3\%} in peak infections (I) over $\mathcal{V}_{\mathcal{N}}$ and a \textbf{3.96\%} reduction in peak hospitalisations (H).\\
Setup: $ c_2^1 \in \{ 9.5,10,10.5,11.0 \} , c_2^2 \in \{ 29.5,30,30.5,31.0 \} , c_2^3 \in \{ 19.5,20,20.5,21.0 \} , c_2^4 \in \{ 24.5,25,25.5,26.0 \} $ and mobility matrix $\text{M}_4$  with $\text{M}_4[1,2] = 0$ and $\text{M}_4[2,1] = 0$.}
\label{fig:rank_13_seqihr}
\end{figure}

\clearpage
\section{Product of Experts (PoE) approach}
\SetVertexStyle[TextFont= \Large]
\begin{figure}[ht]
\vspace{-0.3cm}
\centering
\begin{minipage}{0.46\linewidth}
\resizebox{\linewidth}{!}{%
\begin{tikzpicture}
\begin{pgfonlayer}{background}
 \draw[color =red, style= dashed, rounded corners=0.2cm] (-0.75, -0.6) rectangle (3.3, 2.1) {};
 \draw[color =blue,style= dashed, rounded corners=0.2cm] (3.5, -0.6) rectangle (7.5, 2.1) {};
\end{pgfonlayer}

\Vertex[ size= 0.8,opacity =0, label = $\dot{\textbf{x}}$,   x=1.3, y=0 ]{x_dot_ode}
\Vertex[ size= 0.8,opacity =0, label = $\bm{\theta}$,        x=-0.2,   y=0 ]{theta}
\Vertex[ size= 0.8,opacity =0, label = $\textbf{x}$,         x=2.8,   y=0]{x_ode}
\Vertex[ size= 0.8,opacity =0, label = $\gamma$,             x=1.3, y=1.5 ]{gamma}
\Edge[Direct=true,color= black, lw =1pt](theta)(x_dot_ode)
\Edge[Direct=true,color= black,lw =1pt](gamma)(x_dot_ode)
\Edge[Direct=true,color= black,lw =1pt](x_ode)(x_dot_ode)
\node at (0.0,1.8) {ODE};
\node at (0.0,1.4) {model};

\Vertex[ size= 0.8,opacity =0, label = $\dot{\textbf{x}}$,   x=4, y=1.5 ]{x_dot}
\Vertex[ size= 0.8,opacity =0, label = $\textbf{x}$,         x=5.5, y=1.5]{x}
\Vertex[ size= 0.8,opacity =0, label = $\textbf{y}$,             x=7, y=1.5 ]{y}
\Vertex[ size= 0.8,opacity =0, label = $\phi$,               x=5.5, y=0 ]{phi}
\Vertex[ size= 0.8,opacity =0, label = $\sigma$,             x=7, y=0 ]{sigma}
\Edge[Direct=true,color= black, lw =1pt](phi)(x)
\Edge[Direct=true,color= black,lw =1pt](phi)(x_dot)
\Edge[Direct=true,color= black,lw =1pt](sigma)(y)
\Edge[Direct=true,color= black,lw =1pt](x)(y)
\Edge[Direct=true,color= black,lw =1pt](x)(x_dot)
\node at (4.0,0) {GP};
\node at (4.1,-0.4) {model};

\Edge[Direct=false,style =dashed, color= black, lw =1pt](x_dot_ode)(x_dot)
\Edge[Direct=false,style=dashed, color= black,lw =1pt](x_ode)(x)

\end{tikzpicture}
}
    \caption{ODE model and Gaussian Process model} \label{supp_fig:ode_gp_model}
\end{minipage}\hfill
\begin{minipage}{0.46\linewidth}
\vspace{0.2cm}
\centering
\resizebox{0.65\linewidth}{!}{%
\begin{tikzpicture}

\Vertex[ size= 0.8,opacity =0, label = $\dot{\textbf{x}}$,   x=4, y=1.5 ]{x_dot}
\Vertex[ size= 0.8,opacity =0, label = $\gamma$,             x=4, y=0 ]{gamma}
\Vertex[ size= 0.8,opacity =0, label = $\bm{\theta}$,        x=2.5,  y=1.5 ]{theta}
\Vertex[ size= 0.8,opacity =0, label = $\textbf{x}$,         x=5.5, y=1.5]{x}
\Vertex[ size= 0.8,opacity =0, label = $\textbf{y}$,             x=7, y=1.5 ]{y}
\Vertex[ size= 0.8,opacity =0, label = $\phi$,               x=5.5, y=0 ]{phi}
\Vertex[ size= 0.8,opacity =0, label = $\sigma$,             x=7, y=0 ]{sigma}

\Edge[Direct=true,color= black, lw =1pt](phi)(x)
\Edge[Direct=true,color= black,lw =1pt](phi)(x_dot)

\Edge[Direct=true,color= black,lw =1pt](gamma)(x_dot)
\Edge[Direct=true,color= black,lw =1pt](theta)(x_dot)

\Edge[Direct=true,color= black,lw =1pt](sigma)(y)
\Edge[Direct=true,color= black,lw =1pt](x)(y)
\Edge[Direct=true,color= black,lw =1pt](x)(x_dot)
\end{tikzpicture}
}
  \vspace{0.2cm}
    \caption{Product of Experts (PoE) } \label{supp_fig:poe}
\end{minipage}
\end{figure}

A Gaussian process (GP)  prior is placed on $\textbf{x}_k$, with $\bm{\mu}_k$ and $\bm{\phi}_k$ as the hyper-parameters of this GP with kernel matrix $\textbf{C}_{\phi_{k}}$, we can then write:
\begin{equation}\label{supp_eq:gp}
p(\textbf{x}_k|\bm{\mu}_k,\bm{\phi}_k) = \mathcal{N}(\textbf{x}_k | \bm{\mu}_k, \textbf{C}_{\phi_{k}})
\end{equation}

Using the Gaussian process defined on the state variables $\textbf{x}_k$, we can also make predictions about their derivatives  $\dot{\textbf{x}}_k$ as $\textbf{x}_k$ and $\dot{\textbf{x}}_k$ are jointly Gaussian distributed.  We can write the conditional distribution of the state derivatives as :
\begin{equation}\label{supp_eq:graph_2}
p(\dot{\textbf{x}}_k | \textbf{x}, \bm{\mu}_k, \bm{\phi}_k)  = \mathcal{N}(\dot{\textbf{x}}_k |\textbf{m}_k,\textbf{A}_k ), \; \text{ where}
\end{equation}
$\textbf{m}_k = '\textbf{C}_{\phi_{k}}{\textbf{C}_{\phi_{k}}}^{-1} (\textbf{x}_k - \bm{\mu}_k ); \;\textbf{A}_k = \textbf{C}''_{\phi_{k}} - '\textbf{C}_{\phi_{k}}{\textbf{C}_{\phi_{k}}}^{-1}\textbf{C}'_{\phi_{k}} $

Note that $p(\dot{\textbf{x}}_k | \textbf{x}, \bm{\mu}_k, \bm{\phi}_k)$ corresponds to the second, i.e. GP part of the graphical model in fig. \ref{supp_fig:ode_gp_model}.

Using the functional form of the ODEs and with state specific Gaussian additive noise $\gamma_k$, we can write
\begin{equation}\label{supp_eq:graph_1}
    p(\dot{\textbf{x}}_k| \textbf{X},\bm{\theta},\gamma_k) = \mathcal{N}( \dot{\textbf{x}}_k |\textbf{f}_k(\textbf{X}, \bm{\theta}) , \gamma_k\textbf{I})
\end{equation}

where $\textbf{f}_k(\textbf{X}, \bm{\theta}) = [ \text{f}_k(\text{x}(t_1), \bm{\theta} ),\dots,\text{f}_k(\text{x}(t_N), \bm{\theta} )   ]^T$.
Note that \eqref{supp_eq:graph_1}, corresponds to the ODE part of the graphical model in the fig. \ref{supp_fig:ode_gp_model}.
The two models $p(\dot{\textbf{x}}_k | \textbf{x}_k, \bm{\mu}_k, \bm{\phi}_k)$ in \eqref{supp_eq:graph_2} and $p(\dot{\textbf{x}}| \textbf{X},\bm{\theta},\gamma_k)$ in \eqref{supp_eq:graph_1} are now combined using a product of experts (PoE) approach \cite{prod_of_exp}.

The use of PoE to combine the two models (see fig. \ref{supp_fig:poe}) was originally proposed by \cite{calderhead} and also used by \cite{dondelinger}. 
In the PoE approach the two models are combined to get the distribution of $p(\dot{\textbf{x}}_k| \textbf{X}, \bm{\theta}, \bm{\phi}, \gamma_k)$ as follows:
\begin{alignat}{10}
p(\dot{\textbf{x}}_k| \textbf{X}, \bm{\theta}, \bm{\phi}, \gamma_k ) &\propto p(\dot{\textbf{x}}_k | \textbf{x}_k, \bm{\mu}_k, \bm{\phi}_k)p(\dot{\textbf{x}}| \textbf{X}, \bm{\theta},\gamma_k) \phantom{MMMM} \notag \\
&=\mathcal{N}(\dot{\textbf{x}}_k |\textbf{m}_k,\textbf{A}_k )\mathcal{N}(\dot{\textbf{x}}_k | \textbf{f}_k(\textbf{X}, \bm{\theta}) , \gamma_k\textbf{I}) \label{supp_eq:poe}
\end{alignat}

The joint distribution $p(\dot{\textbf{X}}, \textbf{X}, \bm{\theta}, \bm{\phi},\bm{\gamma})$ can be derived as follows:
\begin{alignat}{10}
p(\dot{\textbf{X}}, \textbf{X}, \bm{\theta}, \bm{\phi},\bm{\gamma}) &= p(\dot{\textbf{X}} | \textbf{X}, \bm{\theta}, \bm{\phi},\bm{\gamma}) p(\textbf{X}| \bm{\phi}) p(\bm{\theta}) p(\bm{\phi}) p(\bm{\gamma})  \notag \\
&= p(\bm{\theta}) p(\bm{\phi})  p(\bm{\gamma}) \prod_{k} p(\dot{\textbf{x}}_k| \textbf{X}, \bm{\theta},\bm{\phi},\gamma_k) p({\textbf{x}}_k| \bm{\phi}_k) \label{supp_eq:joint}
\end{alignat}
where $p(\bm{\theta}), p(\bm{\phi}), p(\bm{\gamma})$ are the prior distributions on ODE parameters $\bm{\theta}$, Kernel hyper-parameters $\bm{\phi}$ and error term $\bm{\gamma}$ respectively. Note that $\bm{\gamma}$ controls the tightness of the gradient coupling.

In \eqref{supp_eq:joint}, after  substituting $p(\dot{\textbf{x}}_k| \textbf{X}, \bm{\theta},\bm{\phi},\gamma_k)$ with \eqref{supp_eq:poe} and $p({\textbf{x}}_k| \bm{\phi}_k)$ with \eqref{supp_eq:gp} we get
\begin{alignat}{10}
p(\dot{\textbf{X}}, \textbf{X}, \bm{\theta}, \bm{\phi},\bm{\gamma}) \propto \;\; & p(\bm{\theta}) p(\bm{\phi}) p(\bm{\gamma}) \times
\prod_{k} \mathcal{N}(\dot{\textbf{x}}_k |\textbf{m}_k,\textbf{A}_k )\mathcal{N}(\dot{\textbf{x}}_k| \textbf{f}_k(\textbf{X}, \bm{\theta}) , \gamma_k\textbf{I})\mathcal{N}(\textbf{x}_k | \bm{\mu}_k, \textbf{C}_{\phi_{k}}) \label{supp_eq:joint_1}
\end{alignat}
In \eqref{supp_eq:joint_1}, marginalization over $\dot{\textbf{X}}$ results in the following:
\begin{adjustwidth}{-0.5cm}{0cm}
\begin{alignat}{10}
\displaystyle p( \textbf{X}, \bm{\theta}, \bm{\phi},\bm{\gamma}) & =  \int p(\dot{\textbf{X}}, \textbf{X}, \bm{\theta}, \bm{\phi},\bm{\gamma}) d\dot{\textbf{X}} \label{supp_eq:margianl_x_dot}\\
& \propto \; p(\bm{\theta}) p(\bm{\phi}) p(\bm{\gamma}) \times \prod_{k} \Big\{ \mathcal{N}(\textbf{x}_k | \bm{\mu}_k, \textbf{C}_{\phi_{k}})\times  
 \int \mathcal{N}(\dot{\textbf{x}}_k |\textbf{m}_k,\textbf{A}_k ) \mathcal{N}(\dot{\textbf{x}}_k| \textbf{f}_k(\textbf{X}, \bm{\theta}) , \gamma_k\textbf{I}) d\dot{\textbf{x}} \Big\} \label{supp_eq:int_def}
\end{alignat}
\end{adjustwidth}
Note that $p( \textbf{X}| \bm{\theta}, \bm{\phi},\bm{\gamma}) =  \displaystyle \prod_k\Big\{ \cdot \Big\}$ and since the 
integral in \eqref{supp_eq:int_def} is analytically tractable, we get:
\begin{alignat}{3}
& p( \textbf{X}, \bm{\theta}, \bm{\phi},\bm{\gamma}) \propto \;
p(\bm{\theta}) p(\bm{\phi}) p(\bm{\gamma}) p( \textbf{X}| \bm{\theta}, \bm{\phi},\bm{\gamma}) \notag\\
&  p(\textbf{X}| \bm{\theta}, \bm{\phi},\bm{\gamma}) \propto \;
\prod_{k} \Big\{ \text{exp}\Big( -\frac{1}{2}\textbf{x}_k^T\textbf{C}_{\phi_k}^{-1}\textbf{x} \Big) \times \text{exp}\Big( -\frac{1}{2}(\textbf{f}_k-\textbf{m}_k)^T(\textbf{A}_k + \gamma_k\textbf{I})^{-1}(\textbf{f}_k-\textbf{m}_k) \Big)\Big\} \notag \\
& \phantom{p(\textbf{X}| \bm{\theta}, \bm{\phi},\bm{\gamma})} \propto \text{exp}\Big[ -\frac{1}{2} \sum_{k} \Big( \textbf{x}_k^T\textbf{C}_{\phi_k}^{-1}\textbf{x} + (\textbf{f}_k-\textbf{m}_k)^T(\textbf{A}_k + \gamma_k\textbf{I})^{-1}(\textbf{f}_k-\textbf{m}_k) \Big) \Big] \label{supp_eq:x_given_theta_phi_gamma}
\end{alignat}
where $\textbf{f}_k$ denotes $\textbf{f}_k(\textbf{X}, \bm{\theta}, \textbf{t})$.
The joint probability for the whole system is given as:
\begin{alignat}{10}
&p(\textbf{Y}, \textbf{X}, \bm{\theta}, \bm{\phi},\bm{\gamma}, \bm{\sigma}) \propto \phantom{}    p(\textbf{Y} | \textbf{X}, \bm{\sigma}) p(\textbf{X}| \bm{\theta}, \bm{\phi},\bm{\gamma})p(\bm{\theta}) p(\bm{\phi}) p(\bm{\gamma}) p(\bm{\sigma}) \label{supp_eq:system_joint}
\end{alignat}
Note that the first factor $ p(\textbf{Y} | \textbf{X}, \bm{\sigma})$ in \eqref{supp_eq:system_joint} is defined in the main text and the second factor $p(\textbf{X}| \bm{\theta}, \bm{\phi},\bm{\gamma})$ is defined in \eqref{supp_eq:x_given_theta_phi_gamma}.
Finally, \cite{dondelinger} developed a Metropolis-Hastings based sampling scheme to sample $\textbf{X}, \bm{\theta}, \bm{\phi},\bm{\gamma}, \bm{\sigma}$ directly from the posterior distribution $p(\textbf{X}, \bm{\theta}, \bm{\phi},\bm{\gamma}, \bm{\sigma}| \textbf{Y}) \propto p(\textbf{Y}, \textbf{X}, \bm{\theta}, \bm{\phi},\bm{\gamma}, \bm{\sigma})$.

In summary, the use of PoE heuristics in literature is common. Both \cite{calderhead} and \cite{dondelinger} using PoE  were able to good estimates for $\bm{\theta}$ using MCMC sampling procedures.

Furthermore, \cite{nico} also using PoE in a Variational inference setting reported a better performance both in terms of computation time and quite surprisingly also in the solution quality, despite their Variational inference setting. The surprising improvement of \cite{nico}  over \cite{calderhead} and \cite{dondelinger} indicated that even though solution approaches developed using PoE as an intermediary step seems to be working empirically, however there is some discrepancy in the overall system. This discrepancy was identified, explained and corrected by \cite{fgpgm}.
\SetVertexStyle[TextFont= \Large]
\begin{figure}[ht]
\centering
\resizebox{0.8\linewidth}{!}{%
\begin{tikzpicture}
\Vertex[ size= 0.8,opacity =0, label = $\dot{\textbf{x}}$,   x=4, y=1.5 ]{x_dot}
\Vertex[ size= 0.8,opacity =0, label = $\gamma$,             x=4, y=0 ]{gamma}
\Vertex[ size= 0.8,opacity =0, label = $\bm{\theta}$,        x=2.5,  y=1.5 ]{theta}
\Vertex[ size= 0.8,opacity =0, label = $\textbf{x}$,         x=5.5, y=1.5]{x}
\Vertex[ size= 0.8,opacity =0, label = $\textbf{y}$,             x=7, y=1.5 ]{y}
\Vertex[ size= 0.8,opacity =0, label = $\phi$,               x=5.5, y=0 ]{phi}
\Vertex[ size= 0.8,opacity =0, label = $\sigma$,             x=7, y=0 ]{sigma}

\Edge[Direct=true,color= black, lw =1pt](phi)(x)
\Edge[Direct=true,color= black,lw =1pt](phi)(x_dot)

\Edge[Direct=true,color= black,lw =1pt](gamma)(x_dot)
\Edge[Direct=true,color= black,lw =1pt](theta)(x_dot)

\Edge[Direct=true,color= black,lw =1pt](sigma)(y)
\Edge[Direct=true,color= black,lw =1pt](x)(y)
\Edge[Direct=true,color= black,lw =1pt](x)(x_dot)

\draw [-> , line width=0.2em ](9.0,0.75) -- (10.5,0.75);
\node[text width=1cm] at (9.8,1.2) {\large{ after}};
\node[text width=4cm] at (10.4,0.1) {\large {marginalization of $\dot{\textbf{x}}$} };

\Vertex[ size= 0.8,opacity =0, label = $\gamma$,             x=14, y=0 ]{gamma1}
\Vertex[ size= 0.8,opacity =0, label = $\bm{\theta}$,        x=12.5,  y=1.5 ]{theta1}
\Vertex[ size= 0.8,opacity =0, label = $\textbf{x}$,         x=15.5, y=1.5]{x1}
\Vertex[ size= 0.8,opacity =0, label = $\textbf{y}$,             x=17, y=1.5 ]{y1}
\Vertex[ size= 0.8,opacity =0, label = $\phi$,               x=15.5, y=0 ]{phi1}
\Vertex[ size= 0.8,opacity =0, label = $\sigma$,             x=17, y=0 ]{sigma1}

\Edge[Direct=true,color= black, lw =1pt](phi1)(x1)


\Edge[Direct=true,color= black,lw =1pt](sigma1)(y1)
\Edge[Direct=true,color= black,lw =1pt](x1)(y1)
\end{tikzpicture}
}
    \caption{Graphical model using PoE (on left) and  the model 
 (on right ) after marginalization of $\dot{\textbf{x}}$. }  \label{supp_fig:poe_marginal}
\end{figure}
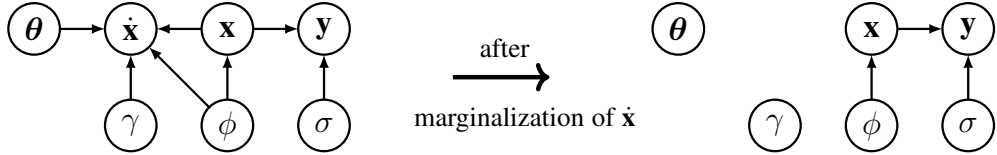

Authors in  \cite{fgpgm} noted that after the marginalization of $\dot{\textbf{X}}$,
the PoE graphical model results in the graphical model shown in \ref{supp_fig:poe_marginal}. In this graphical model there is no link between the parameters $\bm{\theta}$ and the observations $\textbf{y}$, thus defeating the whole purpose of combining the ODE model and GP model (in figure \ref{supp_fig:ode_gp_model}) using PoE. 
Due to this discrepancy, we used the alternative graphical model proposed by \cite{fgpgm} to estimate the posterior-distribution on ODE parameters..

\section{Details on Heuristics}
In Algorithm 1, we use a genetic-algorithm (GA) \cite{nsga2}  as the heuristic to generate the new batch  of policies  (of size $B$)  in step 8, using the objective function values $\{f_{obj}^1,\dots,f_{obj}^B\}$ and constraint violations $\{C_1, \dots,C_B \}$ of the current batch. 

We briefly provide the different steps involved in the working of the GA:
\begin{enumerate}
\item The algorithm starts by randomly generating a set of solutions (also referred as  individuals)  collectively known as the initial population or parent population of size B.
\item The B solutions in the parent population are then evaluated (preferably in parallel), i.e. the objective function values $\{f_{obj}^1,\dots,f_{obj}^B\}$ and constraint violations $\{C_1, \dots,C_B \}$  are computed for each of these individuals. 
\item The individuals  in the parent population are then ranked based on  a criteria which depends on the objective function value and the amount of constraint violation.
\item From the B solutions in the parent population another set of B solutions (or individuals) are generated using genetic operators i.e. SBX crossover and polynomial mutation. The new set of B solutions are collectively known as the child population. 
\item The child population is then evaluated (preferably in parallel) i.e. their objective function values and constraint violations are computed. 
\item The child population and the parent population are then simply merged to form a single population of size 2B known as the mixed population. 
\item The individuals in the mixed population are then assigned a rank based on the ranking criteria. Out of the 2B individual in the mixed population, top B ranking individuals are selected and these selected individuals form the parent population for the next iteration. If the termination criteria (typically a max. number of iterations) is not met, then repeat from step 3.
\end{enumerate}

For exact mathematical details on different operators like SBX crossover, polynominal mutation and the ranking criteria we refer the reader to the original paper \cite{nsga2}. Also, the source code for this GA can be easily obtained from the author's website.

 Following are the values of the different hyper-parameter involved in the GA: population size (B) = 100, probability of cross-over = 0.9, probability of mutation = 0.5, distribution index for  SBX crossover ($\eta_c$) = 10, distribution index for real variable polynomial mutation ($\eta_m$) =  10.
 For the experiments when $\mathbb{K}=3$, we used a maximum of 200 iterations and when $\mathbb{K}=4$, we used a maximum of 250 iterations.

\section{Convergence of Algorithm 1}
Convergence of Algorithm 1 on Nominal $\mathcal{NF}$ and  Stochastic $\mathcal{SF}$ optimization problems for different experiments in the main text. We report mean and standard deviation of 5 runs.
\begin{figure}[h]
\vspace{0.5cm}
\centering
\begin{subfigure}{0.5\linewidth}
    \includegraphics[width = 0.98\linewidth]{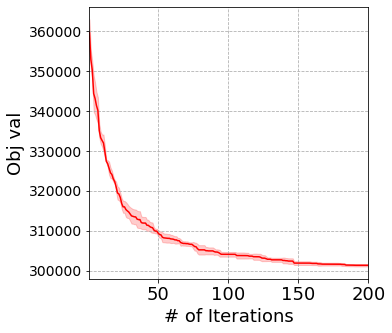}
    \caption{Nominal ($\mathcal{NF}$) problem}    
\end{subfigure}%
\begin{subfigure}{0.5\linewidth}
    \includegraphics[width = 0.98\linewidth]{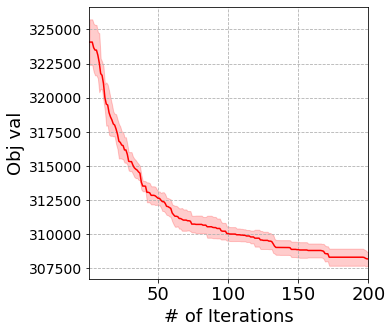}
    \caption{Stochastic  ($\mathcal{SF}$) problem}    
\end{subfigure}
\caption{SEIR model experiment with $\mathbb{K} =3$. }
\end{figure}
\begin{figure}[h]
\centering
\begin{subfigure}{0.5\linewidth}
    \includegraphics[width = 0.98\linewidth]{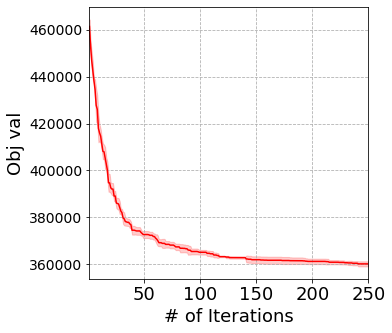}
    \caption{Nominal ($\mathcal{NF}$) problem}    
\end{subfigure}%
\begin{subfigure}{0.5\linewidth}
    \includegraphics[width = 0.98\linewidth]{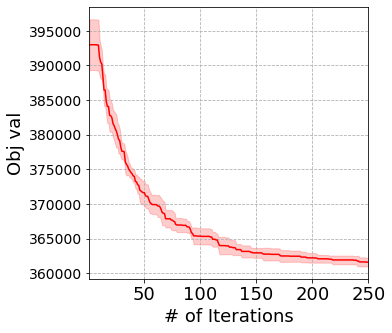}
    \caption{Stochastic  ($\mathcal{SF}$) problem}    
\end{subfigure}
\caption{SEIR model experiment with $\mathbb{K} =4$. }
\end{figure}
\phantom{MM}
\begin{figure}[t!]
\centering
\begin{subfigure}{0.5\linewidth}
    \includegraphics[width = 0.98\linewidth]{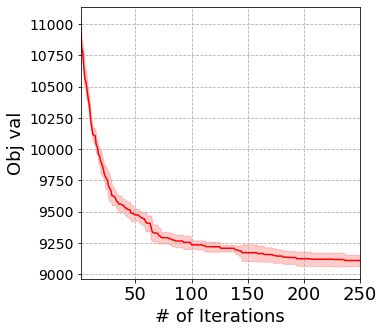}
    \caption{Nominal ($\mathcal{NF}$) problem}    
\end{subfigure}%
\begin{subfigure}{0.5\linewidth}
    \includegraphics[width = 0.98\linewidth]{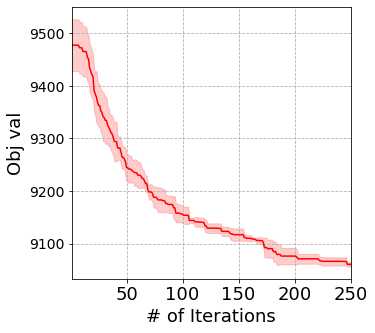}
    \caption{Stochastic  ($\mathcal{SF}$) problem}    
\end{subfigure}
\caption{SEPIHR model experiment with $\mathbb{K} =4$.}
\end{figure}
\newpage
\section{Computational resources}
Our approach consists of the following computational components: i) MCMC sampling, ii) Tractable scenario-set construction with k-means clustering, iii) Solving the nominal vaccine allocation optimization problem and  iv) Solving the stochastic vaccine allocation optimization problem. For the first three components, we run all our experiments on a computer (desktop) with Intel Core i7-6800K CPU with 12 cores and 64GB RAM.  
For the last component, we run all our experiments on a single server node with  Intel(R) Xeon(R) Platinum 8260 CPU with 96 cores and 128 GB RAM.

\clearpage
\section{Tractable scenario-set construction: Effect of changing k}
\subsection{SEIR model}
\begin{figure}[h!]
\vspace{-0.35cm}
\centering
\begin{subfigure}{\textwidth}
  \centering
  \includegraphics[width=\linewidth]{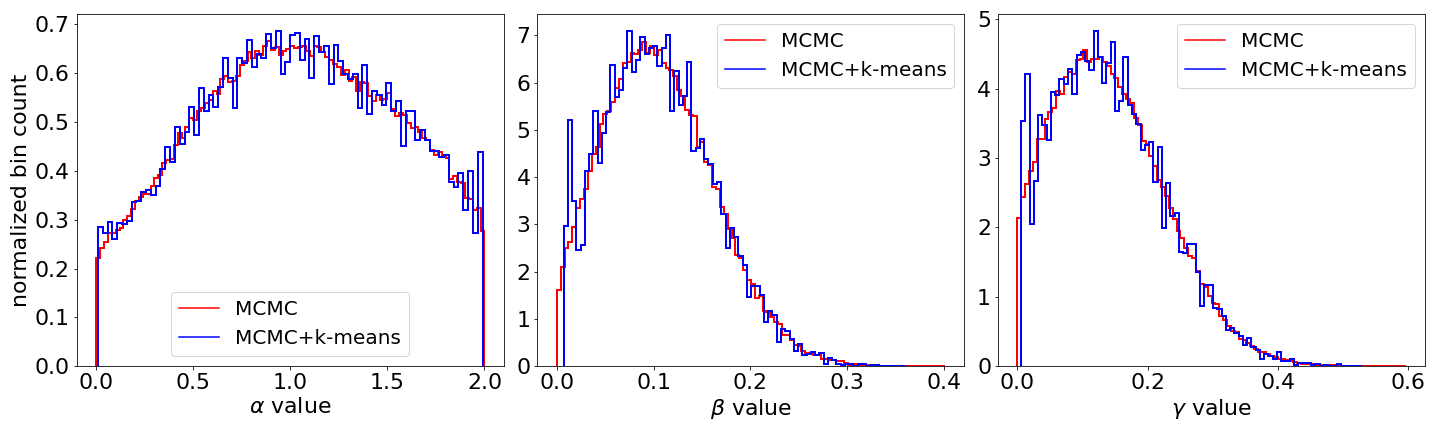}
  \vspace{-0.6cm}
  \caption{k =$0.5\times 10^4$}
\end{subfigure}
\begin{subfigure}{\textwidth}
  \centering
  \includegraphics[width=\linewidth]{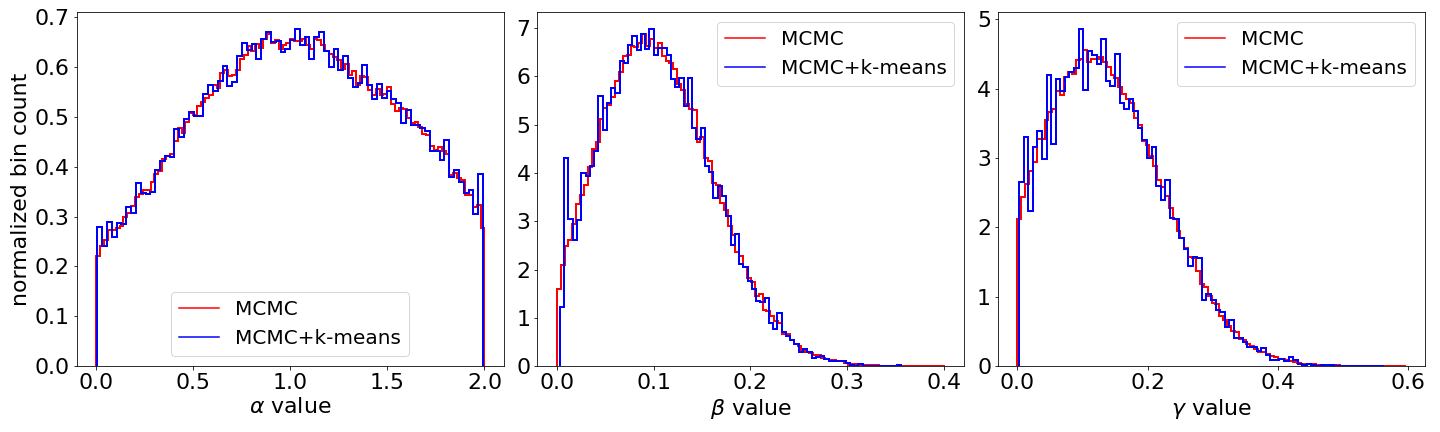}
  \vspace{-0.6cm}
  \caption{k =$1.0\times 10^4$}
\end{subfigure}
\begin{subfigure}{\textwidth}
  \centering
  \includegraphics[width=\linewidth]{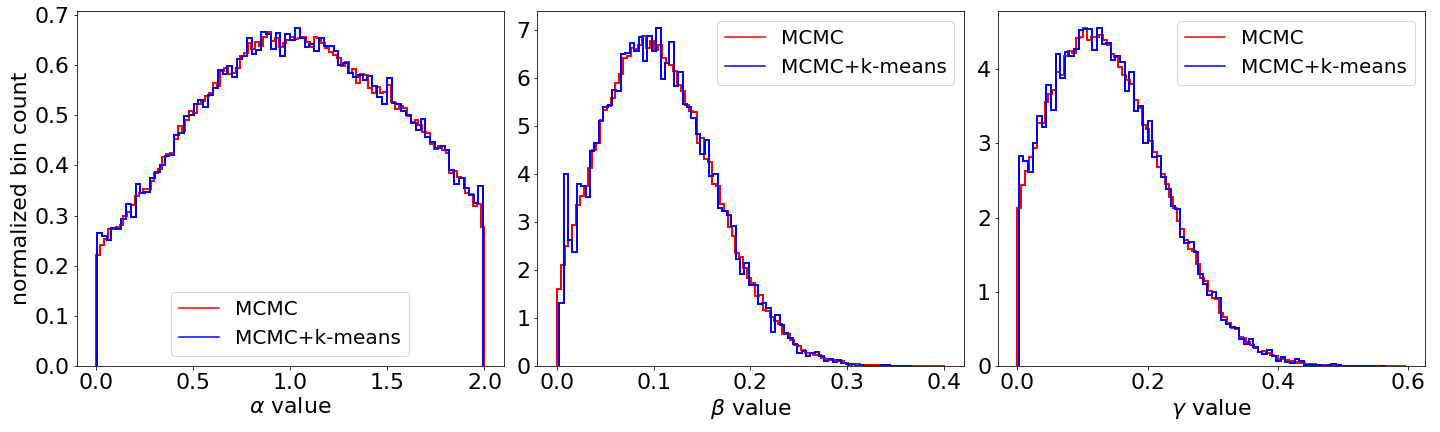}
  \vspace{-0.6cm}
  \caption{k =$1.5\times 10^4$}
\end{subfigure}
\begin{subfigure}{\textwidth}
  \centering
  \includegraphics[width=\linewidth]{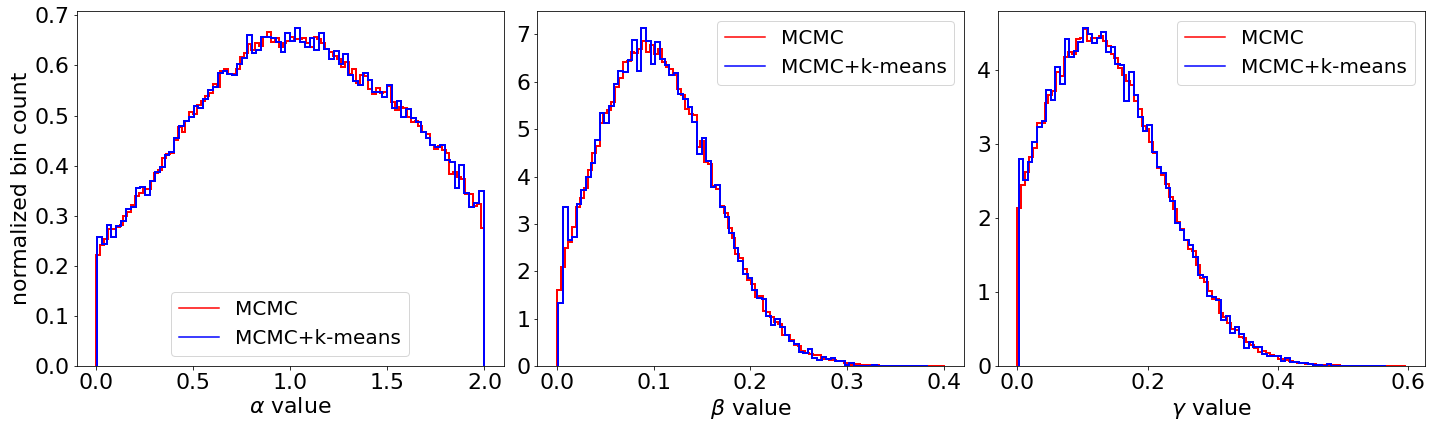}
  \vspace{-0.6cm}
\caption{k =$2.0\times 10^4$}
\end{subfigure}
\end{figure}

\begin{figure}[ht!]\ContinuedFloat
    \centering
\begin{subfigure}{\textwidth}
  \centering
  \includegraphics[width=\linewidth]{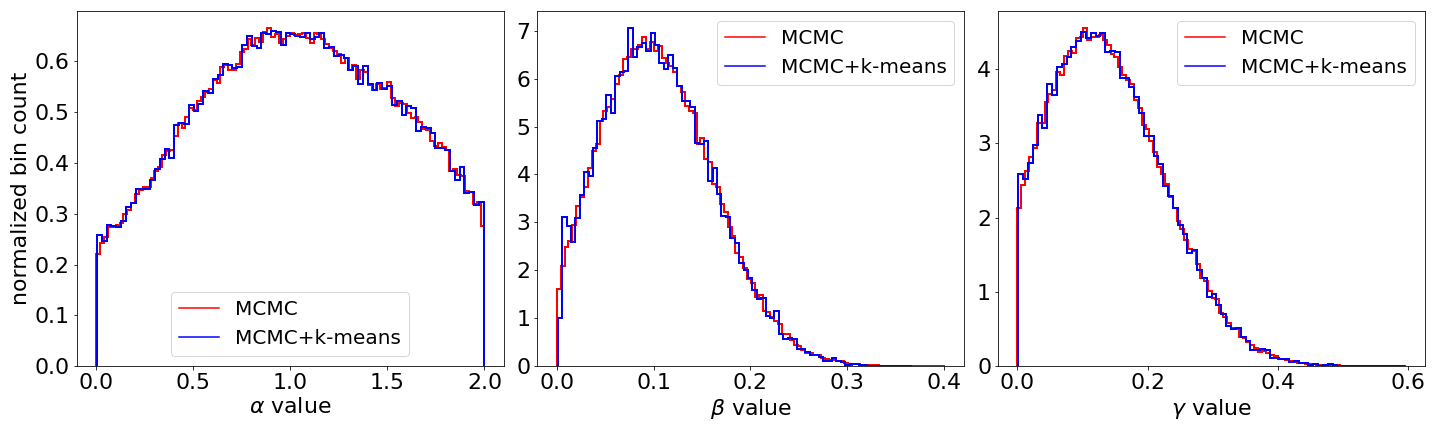}
  \vspace{-0.6cm}
  \caption{k =$2.5\times 10^4$}
\end{subfigure}
\begin{subfigure}{\textwidth}
  \centering
  \includegraphics[width=\linewidth]{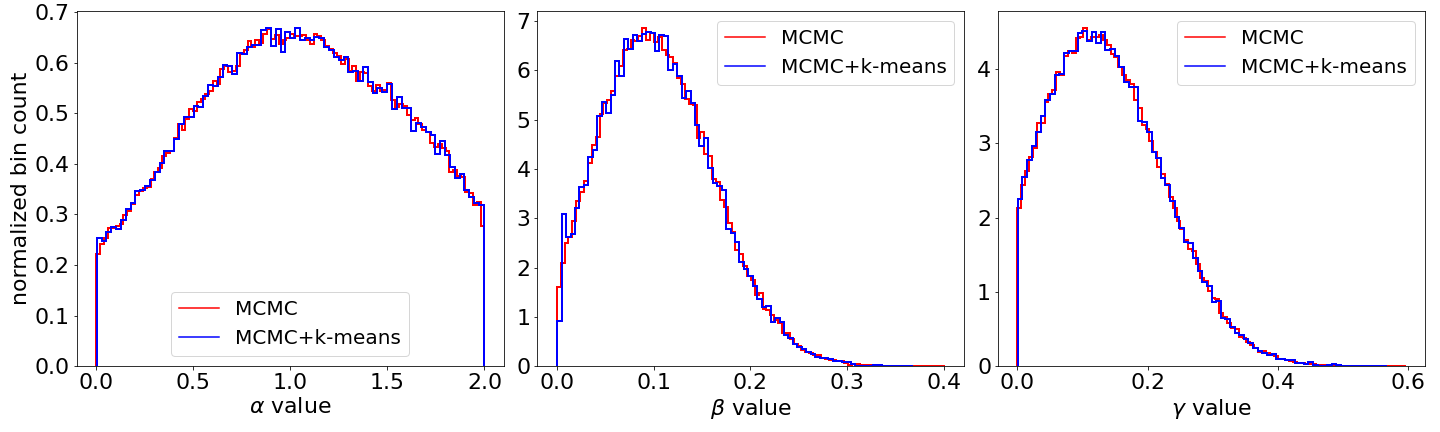}
  \vspace{-0.6cm}
  \caption{k =$3.0\times 10^4$}
\end{subfigure}
\begin{subfigure}{\textwidth}
  \centering
  \includegraphics[width=\linewidth]{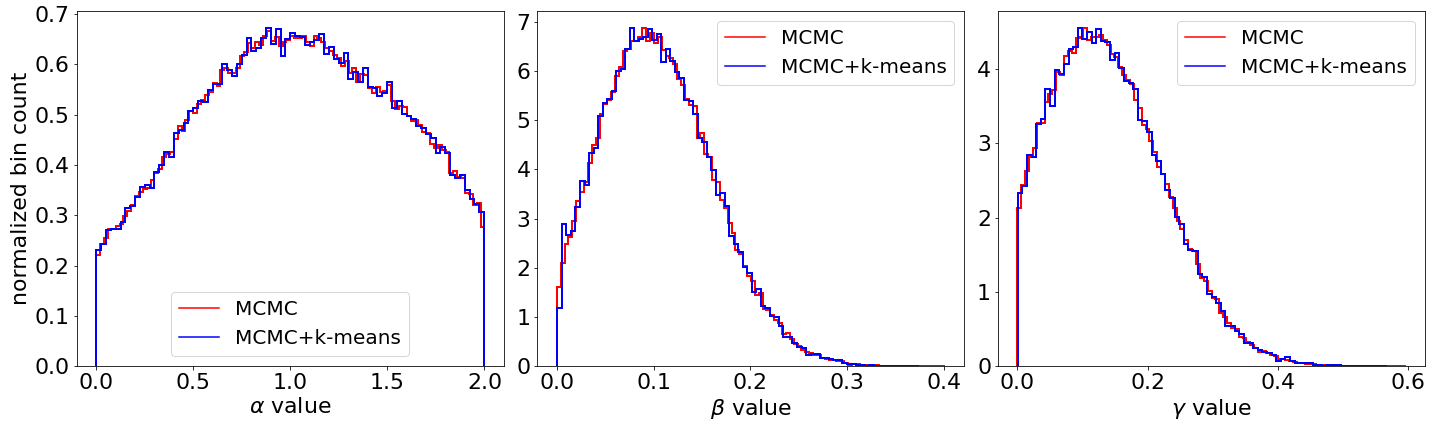}
  \vspace{-0.6cm}
  \caption{k =$3.5\times 10^4$}
\end{subfigure}
\begin{subfigure}{\textwidth}
  \centering
  \includegraphics[width=\linewidth]{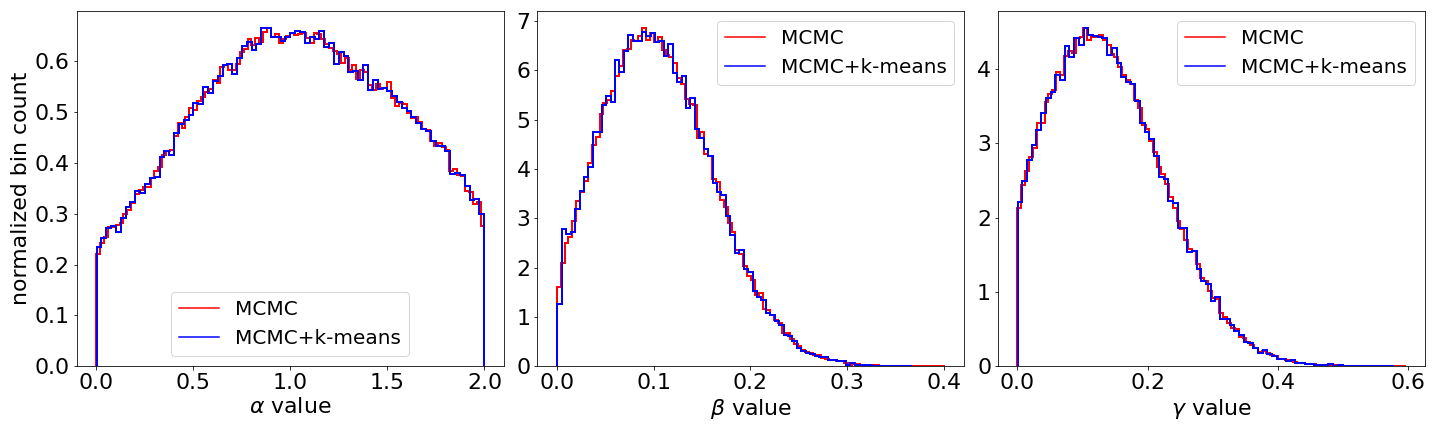}
  \vspace{-0.6cm}
  \caption{k =$4.0\times 10^4$}
\end{subfigure}
\caption{MCMC+k-means denotes the empirical distribution with samples obtained after doing k-means clustering on the original $3 \times 10^5$ samples for the \textbf{SEIR} model. } \label{fig:kmeans_seir_supp}

\end{figure}
\clearpage
\subsection{SEPIHR model}
\centering
\begin{figure}[h!]
\centering
\begin{adjustwidth}{-1.6cm}{-1.6cm}
\begin{subfigure}{\linewidth}
  \centering
  \includegraphics[width=\linewidth]{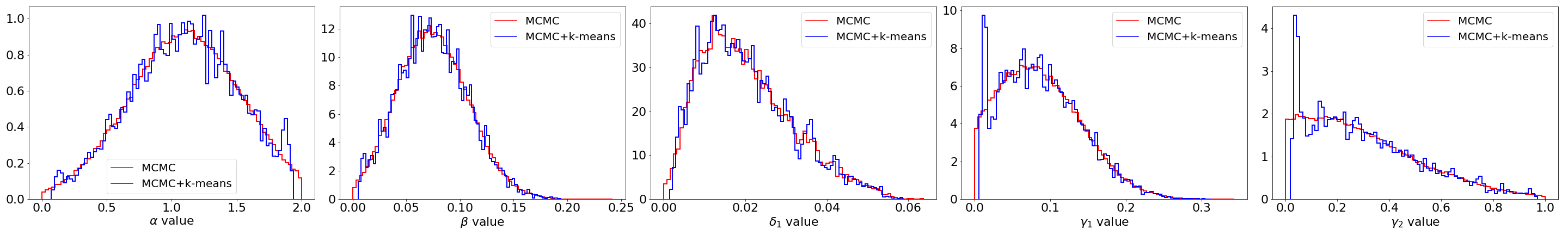}
  \caption{k =$0.5\times 10^4$}
    \vspace{0.5cm}
\end{subfigure} 
\begin{subfigure}{\linewidth}
  \centering
  \includegraphics[width=\linewidth]{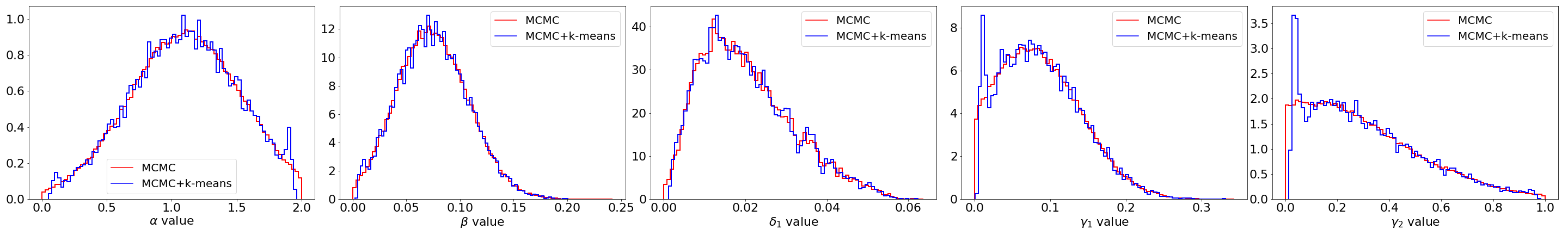}
  \caption{k =$1.0\times 10^4$}
      \vspace{0.5cm}
\end{subfigure}
\begin{subfigure}{\linewidth}
  \centering
  \includegraphics[width=\linewidth]{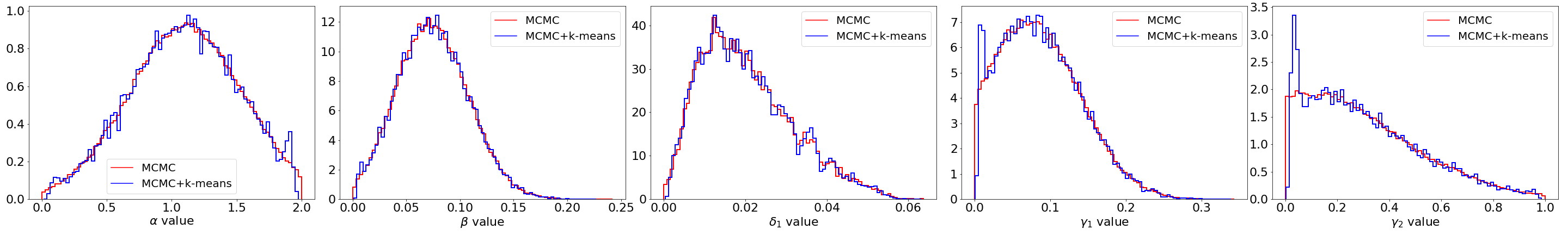}
  \caption{k =$1.5\times 10^4$}
      \vspace{0.5cm}
\end{subfigure}
\begin{subfigure}{\linewidth}
  \centering
  \includegraphics[width=\linewidth]{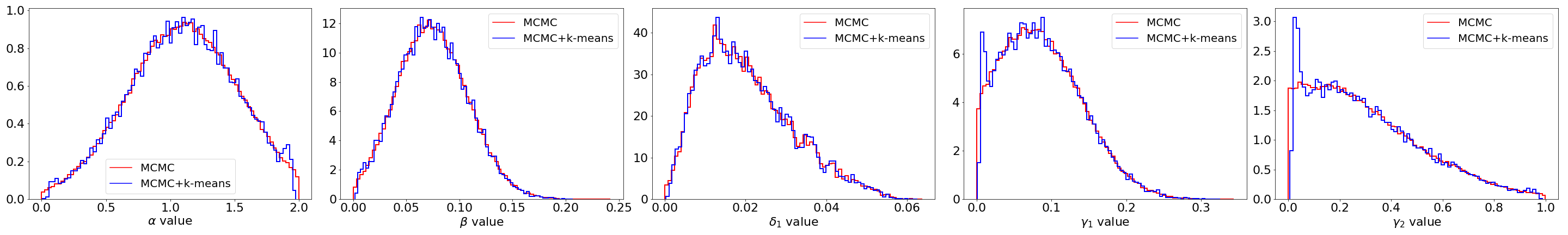}
  \caption{k =$2.0\times 10^4$}
      \vspace{0.5cm}
\end{subfigure}
\begin{subfigure}{\linewidth}
  \centering
  \includegraphics[width=\linewidth]{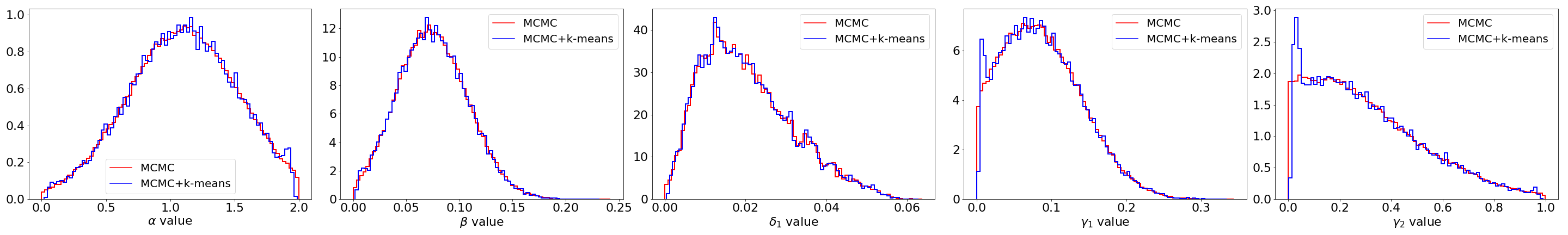}
  \caption{k =$2.5\times 10^4$}
      \vspace{0.5cm}
\end{subfigure}
\end{adjustwidth}
\end{figure}

\begin{figure}[ht!]\ContinuedFloat
\centering
\begin{adjustwidth}{-1.6cm}{-1.6cm}
\begin{subfigure}{\linewidth}
  \centering
  \includegraphics[width=\linewidth]{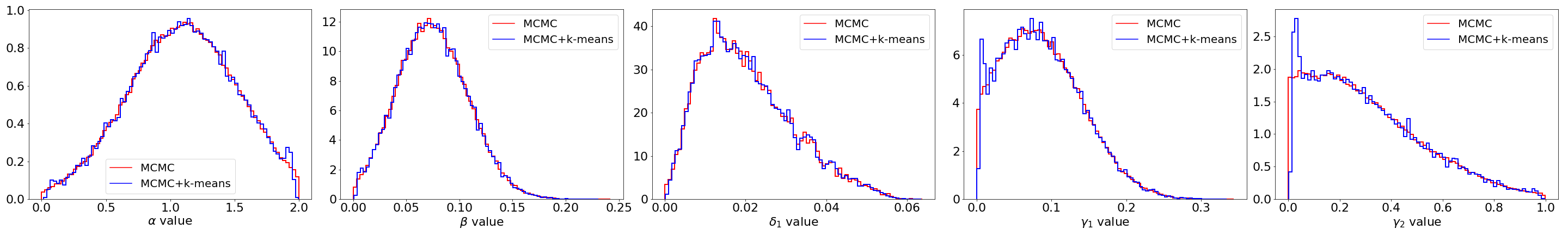}
  \caption{k =$3.0\times 10^4$}
  \vspace{0.5cm}
\end{subfigure}
\begin{subfigure}{\linewidth}
  \centering
  \includegraphics[width=\linewidth]{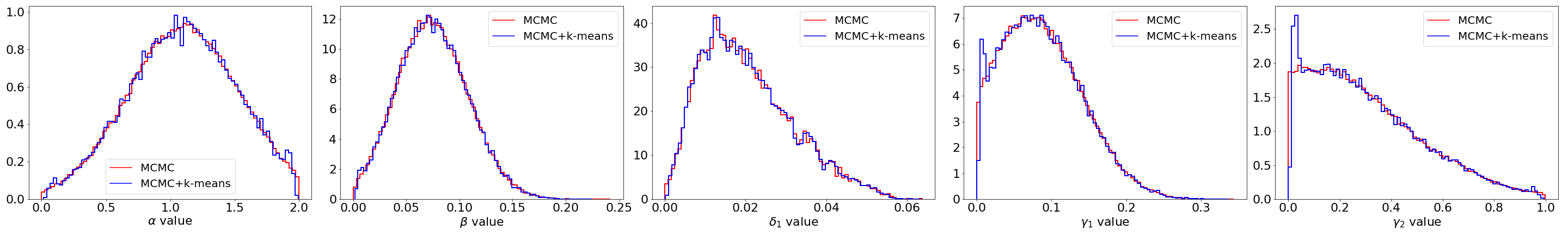}
  \caption{k =$3.5\times 10^4$}
  \vspace{0.5cm}
\end{subfigure}
\begin{subfigure}{\linewidth}
  \centering
  \includegraphics[width=\linewidth]{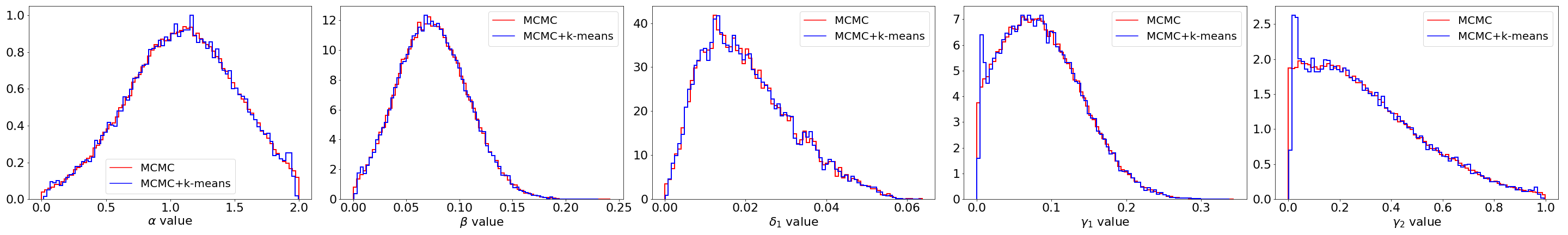}
  \caption{k =$4.0\times 10^4$}
  \vspace{0.5cm}
\end{subfigure}
\begin{subfigure}{\linewidth}
  \centering
  \includegraphics[width=\linewidth]{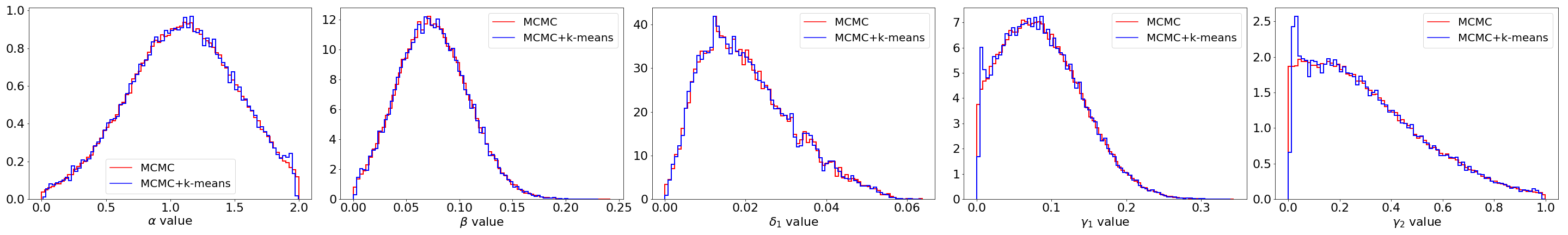}
  \caption{k =$4.5\times 10^4$}
  \vspace{0.5cm}
\end{subfigure}
\begin{subfigure}{\linewidth}
  \centering
  \includegraphics[width=\linewidth]{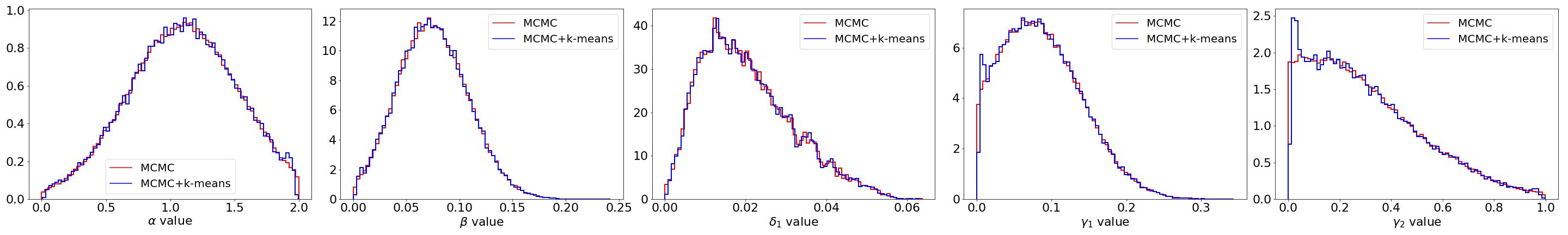}
  \caption{k =$5.0\times 10^4$}
  \vspace{0.5cm}
\end{subfigure}
\caption{MCMC+k-means denotes the empirical distribution with samples obtained after doing k-means clustering on the original $3 \times 10^5$ samples for the \textbf{SEPIHR} model. }
\label{fig:kmeans_sepihr_supp}
\end{adjustwidth}
\end{figure}


\end{document}